\documentclass[a4paper]{article}
\usepackage{imr}
\usepackage{graphicx,bm,amsmath,amsfonts}
\usepackage[english]{babel}
\usepackage{bigints}
\usepackage[dvipsnames]{xcolor}
\usepackage{hyperref}
\usepackage{tikz}
\usepackage{geometry}
\usetikzlibrary{calc}
\usetikzlibrary{patterns}
\usetikzlibrary{arrows}
\usetikzlibrary{arrows.meta}
\usepackage{pgfplots}
\usepackage{placeins}
\usepackage{subfig}
\newcommand{\surface}{\ensuremath{\mathcal{S}}}

\begin{document}

\title{Computing cross fields:\\ A PDE approach based on the Ginzburg-Landau theory}
\author{PA Beaufort$^{1,2}$ \and C Georgiadis$^1$ \and J Lambrechts$^1$\\ \and F Henrotte$^{1,2}$ \and C Geuzaine$^{2}$\and JF Remacle$^1$}
\date{
$^1$UCLouvain -- iMMC, Avenue Georges Lema\^itre 4-6, 1348 Louvain-la-Neuve, Belgium\\
$^2$ULi\`ege -- Montefiore Institute, All\'ee de la D\'ecouverte 10, B-4000 Li\`ege, Belgium 
}

\abstract{
  This paper proposes a method to compute crossfields 
based on the Ginzburg-Landau theory in Magnetism.
According to this theory, 
the magnetic moment distribution in a ferromagnetic material
can be regarded as a vector field with fixed norm, i.e., a directional field.
The energy is the integral over the sample 
of the squared norm of the distribution gradient, 
and the sought distribution is a minimizer 
of this energy under the fixed norm constraint. 
The Ginzburg-Landau functional, which describes mathematically this situation,
has two terms: 
the Dirichlet energy of the distribution
and a term penalizing the mismatch 
between the fixed and actual norm of the distribution. 
Directional fields on surfaces are known to have a number of critical points,
which are properly identified with the Ginzburg-Landau approach:
the asymptotic behavior of Ginzburg-Landau problem provides well-distributed critical points over the 2-manifold, whose indices are as low as possible.
The central idea in this paper is to exploit this theoretical background 
for crossfield computation on arbitrary surfaces. 
Such crossfields are instrumental 
in the generation of meshes with quadrangular elements.
The relation between the topological properties of quadrangular meshes and crossfields
are hence first recalled. 
It is then shown that a crossfield on a surface
can be represented by a complex function of unit norm with a number of critical points,
i.e., a nearly everywhere smooth function 
taking its values in the unit circle of the complex plane. 
As maximal smoothness of the crossfield is equivalent with minimal energy,
the crossfield problem is equivalent 
to an optimization problem based on Ginzburg-Landau functional.
A discretization scheme with Crouzeix-Raviart elements is applied 
and the correctness of the resulting finite element formulation 
is validated on the unit disk by comparison with an analytical solution.
The method is also applied to the 2-sphere where, 
surprisingly but rightly, the computed critical points 
are not located at the vertices of a cube, 
but at those of an anticube.}

\keywords{ Ginzburg-Landau theory, Crossfields, Poincar\'e-Hopf, critical points (singularities), Crouzeix-Raviart, Asterisk Fields}

\maketitle
\thispagestyle{empty}
\pagestyle{empty}

\section{Introduction}

The Finite element method (FEM) provides a powerful and versatile framework for numerical simulation,
which however heavily relies on \emph{mesh generation},
the decomposition of a geometrical region into simple shaped finite elements. 
In two-dimensional geometries, two kinds of elements exist: triangles and \emph{quadrangles}. 
Quadrangular meshes are deemed better than triangular ones because ($i$) there are half as many quadrangles than triangles for the same number of vertices; ($ii$) it is possible to define tensorial operations on quadrangles; and ($iii$) quadrangular meshes ease the tracking of preferred directions in mesh refinement.

However, the generation of quadrangular meshes remains a challenging task, for which many strategies have been explored. Some of them, based on surface parameterization, are suitable for the generation of structured quadrangular meshes, close to regular (square) grids. A crossfield may be used to determine the appropriate parameterization, either on a patch \cite{bommes2009mixed} or globally \cite{kalberer2007quadcover}. A crossfield can also be used for partitioning the surface into a set of curvilinear quadrangular regions (a polyquad), then trivially quadrangulable \cite{kowalski2013pde}. 
The parameterization can also be deduced from a singularity graph \cite{cohen2006designing}. 
In this paper, the primary concern is however to use crossfields as
part of another meshing strategy: a frontal approach firstly proposed
by \cite{lee1994new} that consists in recombining triangles into
quadrangles. This can be done efficiently
\cite{remacle2012blossom} but the quality of the quadrangles
strongly depends on the node location. A heuristic to obtain well
distributed nodes is to spawn them following consistent directions,
such as those suggested by a smooth crossfield. Such a frontal
approach allows building unstructured quadrangular meshes with varying
element size. Other advantages of quadrilateral meshes exist for
specific finite element models: for examples, triangular plate
bending elements are stiffer than quadrilateral ones with the same
number of vertices 

Although there exist various ways to represent discrete crossfields \cite[§5]{vaxman2016directional}, their computation generally relies on some smoothing process, possibly under constraints.
For  an angle-based representation, a crossfield is pictured as four orthogonal or opposite vectors.
From this representation, it is possible to formulate the quadrangulation as a mixed-integer problem \cite{bommes2009mixed}.
More advanced mathematical notions such as holonomy \cite{lai2009metric} may be used as well to design crossfields. This approach requires to build a metric on the 2-manifold. 

In this paper, the so-called Cartesian (complex) representation \cite{palacios2007rotational} is adopted.
This representation naturally takes the symmetries of the cross into account, and the crossfield is identified with a complex-valued function. Complex analysis gives then a large and useful background, 
especially about the theoretical analysis of critical points. 
The second term of the Ginzburg-Landau functional 
is controlled by a parameter depending on the local mesh size. 
When this parameter is small enough, the minimization of the functional results in a smooth crossfield 
whose critical points are optimally located and whose critical points have indices with minimal absolute values, according to the theory.
The previous approach closest to ours is that in \cite{kowalski2013pde}.
It has only the energy term, but the vector field is constrained to have a norm close to the unity.
Critical points are identified in this approach by computing an argument (angle) from the vector field, 
whereas we only need to compute the vector field norm, critical points being in our approach 
points where the crossfield norm locally vanishes. 

Our main contribution is to express the crossfield problem with Ginzburg-Landau equations. Those equations rely on an interesting mathematical and physical backgrounds. In order to grasp the great understanding that Ginzburg-Landau functional provides to the crossfield problem, we first recall the topological constraints of full quadrangular (and triangular) mesh in section \ref{sec:topology} and the link with cross (and asterisk, respectively) field in section \ref{sec:why}. In section \ref{sec:gl}, we develop the intuition of using the Ginzburg-Landau functional for the crossfield problem and we give the related Ginzburg-Landau theory. We derive in section \ref{sec:fem} a simple FEM scheme from the Ginzburg-Landau equations. Our numerical scheme is validated on the unit disk in regards with Ginzburg-Landau theory, section \ref{sec:disk}. On the 2-sphere section \ref{sec:sphere} we get a surprising but correct result. In section \ref{sec:naca}, the Ginzburg-Landau equations are modified to get better results on NACA profiles. Finally, we apply our simple finite scheme on the torus in section \ref{sec:torus}.

\section{Topology of Triangular and Quadrilateral Meshes}\label{sec:topology}

Assume an orientable surface $\surface$ embedded in $\Re^3$.
Let $g$ be the number of handles of the surface.
The topological characteristic $g$, 
which is also called the genus of the surface,
is the maximum number of cuttings along non-intersecting closed curves
that won't make the surface disconnected. 
Let also $b$ be the number of connected components 
of the boundary $\partial\surface$ of the surface.
The Euler characteristic of $\surface$ is then the integer
$$
\chi = 2 - 2g -b.
$$
One has $\chi =2$ for a sphere,
whereas $\chi = 1$ for a disk ($b=1$), 
and $\chi = 0$ for a torus ($g=1$) or a cylinder ($b=2$).

Consider now a mesh on $\surface$ with $n$ nodes (also called vertices), 
$n_e$ edges and $n_f$ facets. 
The Euler formula 
\begin{equation}
\label{eq:chi2}
\chi = n-n_e+n_f
\end{equation} 
provides a general relationship betweeen 
the numbers of nodes, edges and facets in the mesh (details in \cite{eppstein2009nineteen}). 
If $n_b$ nodes (and hence $n_b$ edges) are on the boundary $\partial\surface$,
and if  the number of edges (or nodes) per facet is noted $n_{evf}$
($n_{evf} = 3$ for triangulations and  $n_{evf}=4$ for quadrangulations,
meshes mixing triangles with quadrangles being excluded),
the following identity holds\,:
all facets have $n_{evf}$ edges, 
$n_e-n_b$ edges have two adjacent facets 
and $n_b$ edges have one adjacent facet.
Hence the relationship 
\begin{equation}
\label{eq:rel2}
n_{evf} n_f = 2 (n_e - n_b) + n_b. 
\end{equation}
Elimination of $n_e$  between \eqref{eq:rel2} and \eqref{eq:chi2} yields 
\begin{equation}
\label{eq:rel3}
2n - n_b + (2-n_{evf}) n_f = 2 \chi,
\end{equation}
which is true for any triangulation or quadrangulation.

A regular mesh has only regular vertices.
An internal vertex is regular 
if it has exactly $6$ adjacent triangles or $4$ adjacent quadrangles,
whereas a boundary vertex is regular 
if it has exactly $3$ adjacent triangles or $2$ adjacent quadrangles. 
One has then
\begin{equation}
\label{eq:rel4}
6(n-n_b)+3n_b=3n_f 
\quad \Rightarrow \quad
n_f = 2 n - n_b
\end{equation}
and
\begin{equation}
\label{eq:rel5}
4(n-n_b)+2n_b=4n_f 
\quad \Rightarrow
\quad 
n_f = n - {n_b \over 2}
\end{equation}
respectively for a regular triangulation and a regular quadrangulation, from a topological point of view.
Substitution of \eqref{eq:rel4} and \eqref{eq:rel5} into \eqref{eq:rel3} 
shows that only surfaces with a zero Euler characteristic can be paved with a regular mesh.
If $\chi \neq 0$, irregular vertices will necessarily be present in the mesh.

The number and the index of the irregular vertices
is tightly linked to the Euler characteristic $\chi$,
which is a topological invariant of the surface.
We call valence of a vertex the number of facets adjacent to the vertex in the mesh. 
In a regular mesh, all vertices have the same valence $v_{reg}$.
In a non regular mesh, on the other hand,
a number of irregular vertices have a valence $v \ne v_{reg}$,
and one notes the integer $k = v_{reg} -v$ the valence mismatch of a vertex. 

Assume a quadrangulation with 
$n_k$ irregular internal vertices of valence $v=4-k$,
and $n_{bk}$ irregular boundary vertices of valence $2-k$, $k$ given.
All other vertices are regular.
There are then $n-n_b-n_k$ regular internal vertices of valence $4$,
and $n_b-n_{bk}$ regular boundary vertices of valence $2$,
so that one can write
\begin{equation}
\label{eq:rel6}
4 n_f = 4 (n-n_b-n_k) + 2(n_b-n_{bk}) + (4-k) n_k + (2-k) n_{bk},
\end{equation}
and the substraction of \eqref{eq:rel3} with $n_{evf}=4$  yields
$$
\chi = {k \over 4} (n_k + n_{bk}),
$$
showing that, in a quadrangulation, 
each irregular vertex counts for 
$\text{index}  ( {\mathbf x}_i ) = k/4$ in the Euler characteristic,
a quantity called the indice of the irregular vertex $ {\mathbf x}_i $.

Summing up now on different possible values for $k$,
one can establish that a  quadrangulation 
of a surface with Euler characteristic $\chi$  verifies
\begin{equation}
\label{eq:hopfq}
\chi = \sum_k  {k \over 4} (n_k + n_{bk}) = 
\sum_{i=1}^N \text{index}  ( {\mathbf x}_i ).
\end{equation}

\begin{figure}
\begin{center}
\includegraphics[width=4cm]{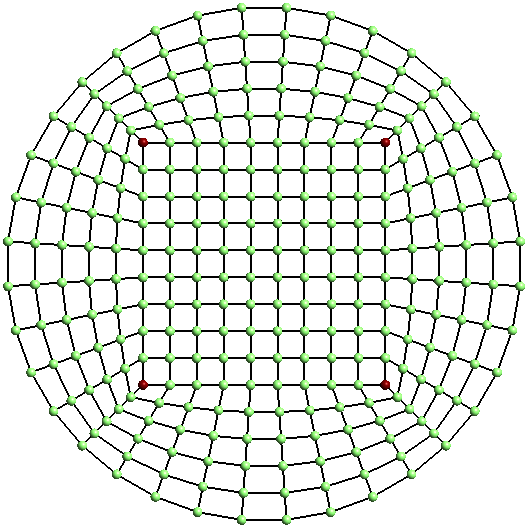}
\includegraphics[width=4cm]{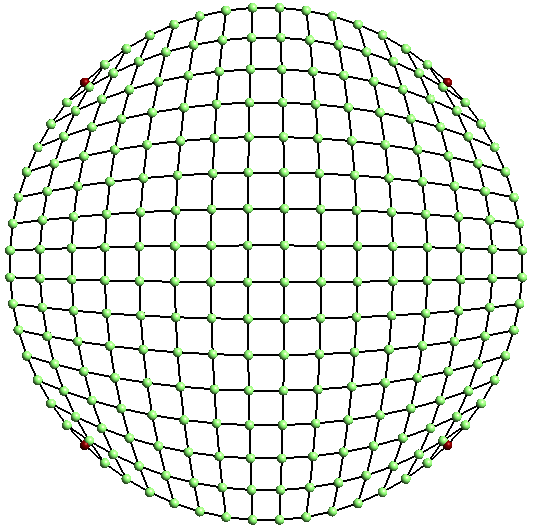}
\end{center}
\caption{A quadrilateral mesh of a circle. Four irregular vertices of index $1/4$ (in red) are required to obtain such a mesh. The irregular vertices may be inside
the disk (left) or on its boundary (right) \label{fig:circle_quads}}
\end{figure}

Consider, for instance, the quadrangulation of a disk,
which is a surface with $\chi=1$. 
A minimum of $n_1=4$ irregular vertices of  index $1/4$
must be present.
They can be located either on the boundary (vertices of valence 1)
or inside the disk (vertices of valence 3), 
Fig.~\ref{fig:circle_quads}. 

\begin{figure}
\begin{center}
\includegraphics[width=2.5cm]{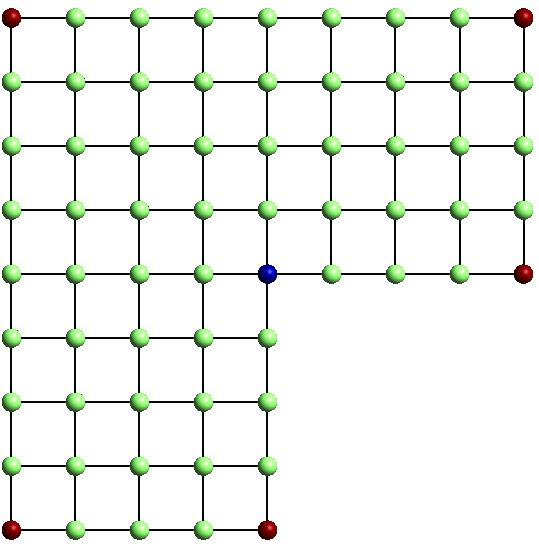}
\includegraphics[width=2.5cm]{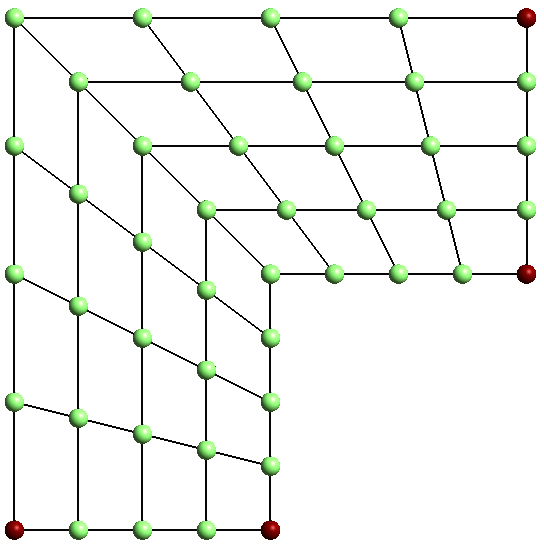}
\includegraphics[width=2.5cm]{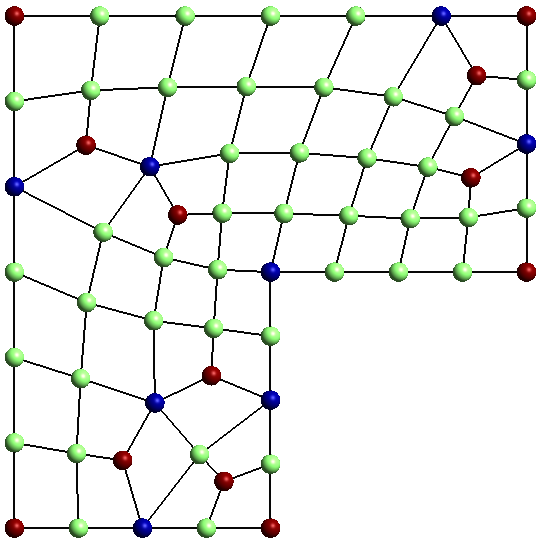}
\end{center}
\caption{Different quadrangulations of a L-shaped domain. 
Irregular vertices of index $1/4$ are displayed in red, 
whereas ones of index $-1/4$ are displayed in blue. 
The sum of the indices of the irregular vertices is equal to $\chi=1$ in all cases.
\label{fig:nonsmooth}}
\end{figure}
Fig.\,\ref{fig:nonsmooth} shows three different quadangulations 
of a L-shaped domain ($\chi=1$).
Regular boundary nodes should all have a valence of 2.
The mesh on the left has $6$ irregular vertices located at the corners of the domain\,:
five with index $1/4$, and one with index $-1/4$. 
The central mesh, on the other hand, 
has the minimum amount of irregular vertices, i.e. four ones of index $1/4$. 
The right mesh generated 
by recombination of a standard Delaunay triangular mesh (see \cite{remacle2012blossom})
has twelve vertices of index $1/4$, and eight vertices of index $-1/4$,
both on the boundary and inside the domain.  
Quality meshes should have as few irregular vertices as possible. 
In what follows, a general approach allowing to compute the position of such irregular vertices 
before meshing the surface is presented.

\section{Why Crossfields?}\label{sec:why}

Crossfields are auxiliary in the generation of quadrangular meshes. 
We shall show that nonregular vertices defined in the previous section
are precisely the critical points of a crossfield,
and that these critical points of the crossfield can also be related 
to the Euler characteristic of the meshed surface.
This result represents an important theoretical limit
on the regularity of quadrangular meshes.

\subsection{Continuity}

A crossfield $f$ is a field defined on a surface $\surface$
with values in the quotient space $S^1/Q$, 
where $S^1$ is the circle group 
and $Q$ is the group of quadrilateral symmetry. 
Pictorially, it associates to each point of the surface $\surface$, which has to be meshed,
a cross made of four unit vectors that are
orthogonal with each others in the tangent plane $T\surface$ of the surface.

\begin{figure}
\begin{center}
\includegraphics{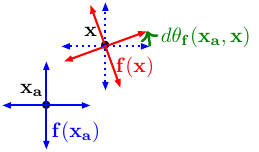}
\end{center}
\caption{Differential function $\mathrm{d}\theta_f$.}
\label{fig:dtheta}
\end{figure}

A surface $\surface$ can be identified with its tangent space 
in any neighborhood $\sigma \subset\surface$
that is sufficiently small to have curvature effects negligible. 
This local identification of the surface with a vector space
endows it with a natural parallel transport rule,
so that 
the angular differential $\mathrm d\theta_f({\mathbf x}_a, {\mathbf x}_b)$ 
can be defined as the minimal angle, with its sign, 
between the branches of $f({\mathbf x}_a)$ 
and any of the branches of $f({\mathbf x}_b)$
for any pair of points ${\mathbf x}_a, {\mathbf x}_b \in\sigma$
where $f$ is defined, Fig.~\ref{fig:dtheta}.
Taking now as reference the cross $f({\mathbf x}_a)$,
an angular coordinate
\begin{equation}
\theta_f({\mathbf x}) = \mathrm d\theta_f({\mathbf x}_a, {\mathbf x})
\label{eq:thetaf}
\end{equation}
can be defined for crosses in $\sigma$.
The crossfield $f$ is deemed continuous (regular) at ${\mathbf x}_b$ 
if the limit
\begin{equation}
\lim_{{\mathbf x} \rightarrow {\mathbf x}_b}  \theta_f({\mathbf x}) = \theta_f({\mathbf x}_b) 
\label{eq:continuity}
\end{equation}
exists (i.e. is unique).
It is then equal to $\theta_f({\mathbf x}_b) $.
Isolated points ${\mathbf x}_i,~i=1\dots N$, of $\surface$
where the limit (\ref{eq:continuity}) does not exist 
are called critical points or zeros of the crossfield. 

\subsection{Index and degree}

Although defined locally, 
the notion of continuity gives unexpectedly 
valuable information about the topology of $\surface$,
which is a nonlocal concept.
To see this, consider a crossfield $f$
defined on a quadrangular element
delimited by four (possibly curvilinear) edges. 
Assume the crossfield is parallel to the four edges 
(i.e. one of the four branches of the cross 
is parallel to the tangent vector of the edge
at each point of the edge, except the extremities)
and prolongates smoothly inside the quadrangle.
This field is discontinuous 
at corners where edges do not meet at right angle,
but it is continuous everywhere else.
Making the same construction for all elements of a quadrangular mesh,
one obtains a crossfield  $f$ 
topologically identified with the quadrangular mesh,
and that is continuous everywhere except at the vertices of the mesh.
This field has thus got isolated critical points at mesh vertices,
but not all critical points have the same significance.
Some critical points have a specific topological value,
associated with the notion of index.

To introduce the notion of index,
an angular coordinate needs to be defined 
for points in a neighborhood $\sigma_i$ of a critical point ${\mathbf x}_i$.
Picking up an arbitrary regular point ${\mathbf x}_a \in \sigma_i$,
${\mathbf x}_a \ne {\mathbf x}_i$,
the local unit vector basis
$$
{\mathbf e}_1 = \frac{{\mathbf x}_a - {\mathbf x}_i}{|{\mathbf x}_a - {\mathbf x}_i|}
\quad , \quad
{\mathbf e}_2 = {\mathbf n} \times {\mathbf e}_1,
$$ 
with $ {\mathbf n}$ the normal to $\surface$,
is constructed, and 
hence a local polar coordinate system
\begin{equation}
r({\mathbf x}) = |{\mathbf x} - {\mathbf x}_i|
\quad , \quad 
\theta({\mathbf x}) = \mathrm{atan2}\Big(
({\mathbf x} - {\mathbf x}_i) \cdot {\mathbf e}_2, 
({\mathbf x} - {\mathbf x}_i) \cdot {\mathbf e}_1 \Big)
\label{eq:theta}
\end{equation}
can be defined for points in $\sigma_i$.

\begin{figure}
\begin{center}
\includegraphics{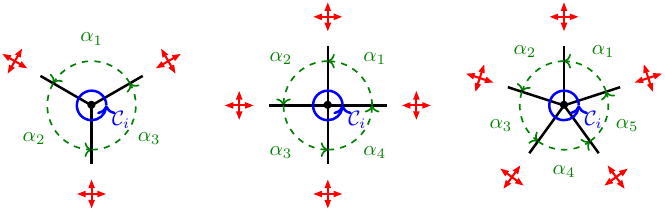}
\end{center}
\caption{Illustration of vertices where the indices of the crossfield (in red) are respectively  $1/4, 0 \text{ and } -1/4$, from left to right.
The index only depends on the number of quadrangles that are adjacent to the vertex,
independently of the values of the angles $\alpha_p$,
which don't need to be identical as they are in the figure.} 
\label{fig:index}
\end{figure}

A circular curve $\mathcal{C}_i$ 
of infinitesimal radius centered around the vertex ${\mathbf x}_i$
is now considered.
As the angles $\theta({\mathbf x})$ (\ref{eq:theta})
and $\theta_f({\mathbf x})$ (\ref{eq:thetaf})
are precisely the elements of the groups $S^1$ and $S^1/Q$,
respectively,
the crossfield on $\mathcal{C}_i$ can be regarded as a mapping
\begin{equation}
f \,:S^1 \mapsto S^1/Q.
\label{eq:s1s1q}
\end{equation}
The mapping is continuous, since
$\mathcal{C}_i$ circles around the critical point ${\mathbf x}_i$,
but it does not cross it.
The index of $f$ at ${\mathbf x}_i$ is the degree of the mapping (\ref{eq:s1s1q}),
i.e. the number of times the codomain  wraps around the domain 
under the mapping.
Its algebraic expression is easily expressed 
in terms of the angles $\theta$ and $\theta_f$ as 
$$
\mathrm{index}({\mathbf x}_i) = 
\frac 1 {2\pi} {\oint_{\mathcal{C}_i}  \mathrm d\theta_f}
$$
where $2\pi$ is $\oint_{\mathcal{C}_i}  \mathrm d\theta$. 
In case of a vertex ${\mathbf x}_i$ of valence $v_i$,
i.e. a vertex adjacent to $v_i$ quadrangular elements,
the integral evaluates as
\begin{equation}
\mathrm{index}({\mathbf x}_i) 
= \frac 1 {2\pi} \sum_{p=1}^{v_i}  \left(\alpha_p - \frac \pi 2\right) 
= \frac 1 {2\pi}  (2\pi-v_i\frac \pi 2) = \frac {4-v_i} 4,
\label{eq:dtheta}
\end{equation}
where the $\alpha_p$'s 
are the angles of the $v_i$ quadrangular elements
adjacent to the considered vertex ${\mathbf x}_i$,
and where the obvious relationship 
$\sum_{p=1}^{v_i} \alpha_p = 2\pi$ has been used.
The crossfield $f$ has index 0 
at vertices adjacent to four quadrangular elements, 
whereas it has index $1/4$ ($-1/4$)
at vertices adjacent to 3  (5, respectively) quadrangular elements meet, Fig.~\ref{fig:index}.
As one sees, the index
is a topological characteristic of the crossfield $f$ 
at the critical point ${\mathbf x}_i$.
It does not depend on the choice of the  curve $\mathcal{C}_i$,
nor on the choice of an angular reference
for the angles $\theta({\mathbf x})$ and $\theta_f({\mathbf x})$.

\subsection{Poincar\'e-Hopf theorem}

Equation (\ref{eq:dtheta})
relates the index of the crossfield at a critical point ${\mathbf x}_i$
with one fourth of valence $k_i=4-v_i$ of the corresponding mesh vertex. 
This result can be combined with the algebraic topology result 
of previous section (\ref{eq:hopfq}) 
that each internal irregular vertex of valence $k_i$ counts for $k_i/4$ 
in the Euler characteristic of the underlying surface. 
This yields the relationship
\begin{equation}
\sum_{i=1}^N \text{index}  ( {\mathbf x}_i ) = \chi
\label{eq:poincarehopf}
\end{equation}
for the critical points of a crossfield $f$ defined on a surface $\surface$. 

This is a generalization Poincar\'e-Hopf theorem,
which states that the sum of the indices of the critical points
of a {\em vector} field ${\mathbf v}$ defined on a surface $\surface$
without boundary is equal to the Euler characteristic of the surface. 
This famous theorem draws an unexpected and profound link 
between two apparently distinct areas of mathematics,
topology and analysis. 
Whereas vector fields have integer indices at critical points,
crossfields have indices that are multiples of 1/4. 
Still the topological relationship (\ref{eq:poincarehopf})
of Poincar\'e-Hopf holds in both cases. Actually, our developments reach same inferences as \cite{ray2008n}.

\section{Crossfield Computation: the Planar Case}\label{sec:gl}

We introduce the representation of a crossfield by means of a vector field.
From this representation, we derive the problem to solve that corresponds to minimize Ginzburg-Landau functional.
Its asymptotic behavior provides suitable critical points, if any.

\subsection{Vector representation of crossfields}

Only scalar quantities can be compared at different points of a manifold.
For the comparison or, more generally, 
for differential calculus with nonscalar quantities like crossfields,
a parallel transport rule needs to be defined on the manifold.  
On a surface (two-manifold), this rule can take the form
of a regular vector field
which gives at each point the direction of the reference angle 0.
Poincar\'e-Hopf theorem says 
that such a field does not exist in general,
and in particular on manifolds whose Euler characteristic is not zero.
The situation is however easier in the planar case.
A global Cartesian coordinate frame can always be defined over the plane,
and be used to evaluate the orientation of the crossfield. 
We shall therefore expose the crossfield computation method
in the planar case,
and then generalize to nonplanar surfaces,
where we will have to deal with local reference frames, 
in a subsequent section.

A cross $f({\mathbf x})$ is an element of the group $S^1/Q$,
which can be represented by the angle $\theta_f({\mathbf x})$ 
it forms with the local reference frame. 
Yet, due to the quadrilateral symmetry, 
four different angles in $[0,2 \pi[$ represent the same crossfield $f({\mathbf x})$. 
Let for instance the angles $\theta_1 = 0$ and $\theta_2=\pi/2$ 
represent the same cross.
The average $(\theta_1+\theta_2)/2 = \pi/4$ represents another cross,
whereas the difference $\theta_2-\theta_1 = \pi/2$ is not zero.
So, we have $\frac 1 2 (x + x) \ne x$ and $x-x \ne 0$,
which clearly indicates that the values of the crossfield $f$
do not live in a linear (affine) space.
This makes the representation by $\theta_f$ 
improper for finite element interpolation. 
The solution is two-fold.
First, the angle $\theta_f$ is multiplied by four,
so that the group $S^1/Q$ is mapped on the unit circle $S^1$,
and the cross $f$ is therefore represented by a unit norm vector ${\mathbf f}$.
Then, the vector is represented in components
in the reference frame as 
$$
{\mathbf f} = (\cos 4\theta_f, \sin 4 \theta_f) \equiv (f_1, f_2).
$$

This vector may be represented by a complex-valued function
$$
f = f_1 + i ~ f_2
$$
This representation corresponds to a vector field that is described by a complex exponential whose argument is $4~\theta$.
A crossfield is thus depicted by the fourth roots of a (unit) complex number.
This observation may be generalized for directional fields with $n$ symmetries \cite[§5.2]{vaxman2016directional}.

\subsection{Laplacian smoothing}

Computing the crossfield $f$ consists thus now 
of computing the vector field representation ${\mathbf f}$,
which obviously lives in a linear space (a 2D plane). 
The components of ${\mathbf f}$ are fixed 
on the boundaries of $\Gamma = \partial \surface$
so that the crosses are parallel with the exterior normal vector 
${\mathbf n} = (\cos \theta_n, \sin \theta_n)$
i.e.
$$
{\mathbf f} = (\cos 4\theta_n, \sin 4 \theta_n)
\quad \mathrm{ on } \ \Gamma.
$$
Propagating ${\mathbf f}$ inside $\surface$ is here done by solving a Laplacian problem. 
Even though the vector representation ${\mathbf f}$ is unitary on $\Gamma$, 
it tends to drift away from $S^1$ inside the domain.
The computed finite element solution ${\mathbf f}$ 
lies therefore outside the unit circle and must be projected back on $S^1$
to recover the angle
$$
\theta_f = {\mathrm{atan2} (f_2, f_1) \over 4}.
$$

Due to the multiplication by 4,
the indices of the critical points of the vector field ${\mathbf f}$ verify
\begin{equation}
\sum_{i=1}^N \text{index}  ( {\mathbf x}_i ) = 4 \chi.
\label{eq:phtimes4}
\end{equation}

\subsection{The Ginzburg-Landau model}\label{sec:glm}

\begin{figure}[ht!]
\begin{center}
\subfloat[Minimizing Dirichlet energy.]{\includegraphics[scale=.1]{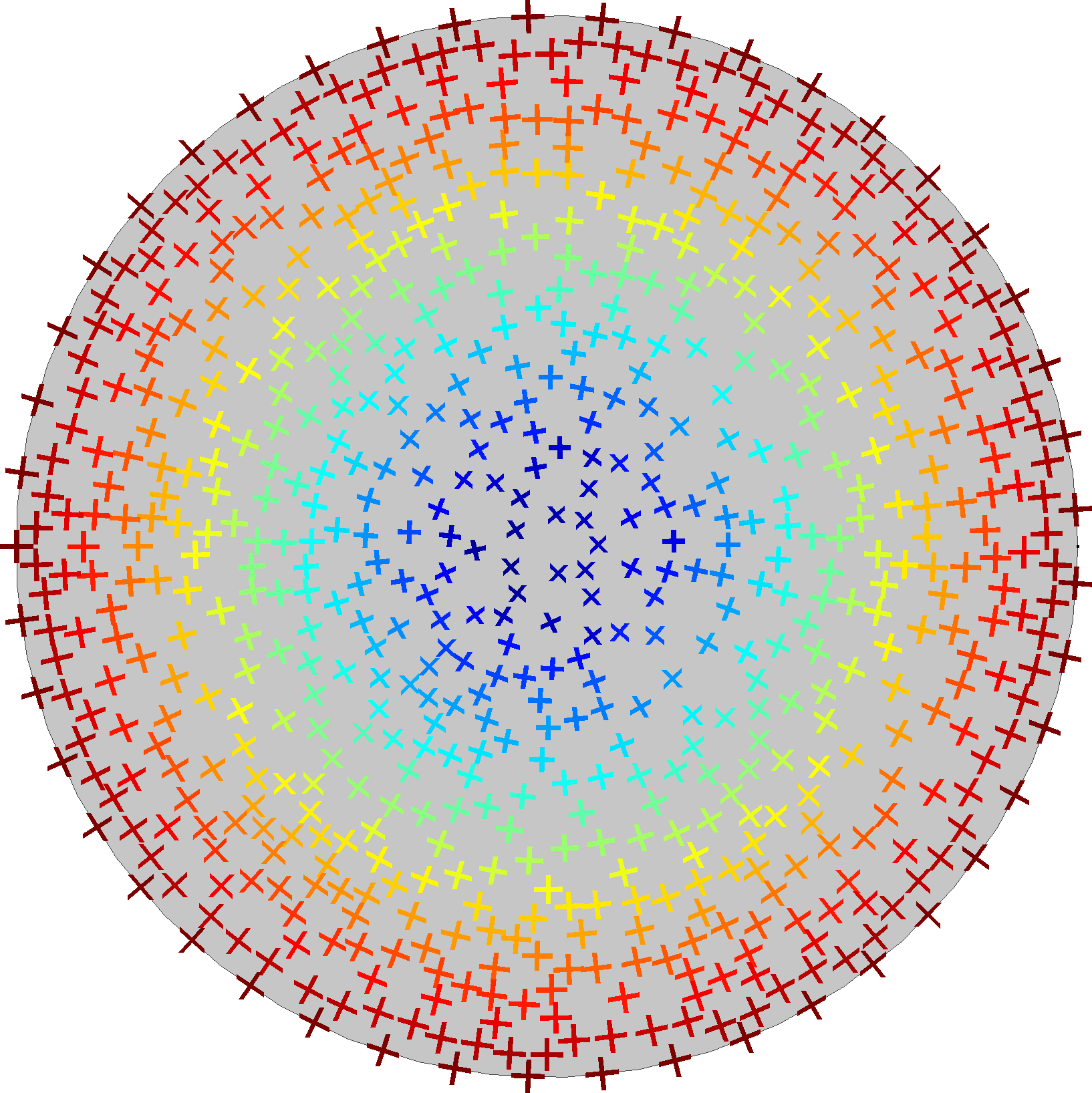}\label{sub:diskdir}}\hspace{.25cm}
\subfloat[Minimizing Ginzburg-Landau energy.]{\includegraphics[scale=.1]{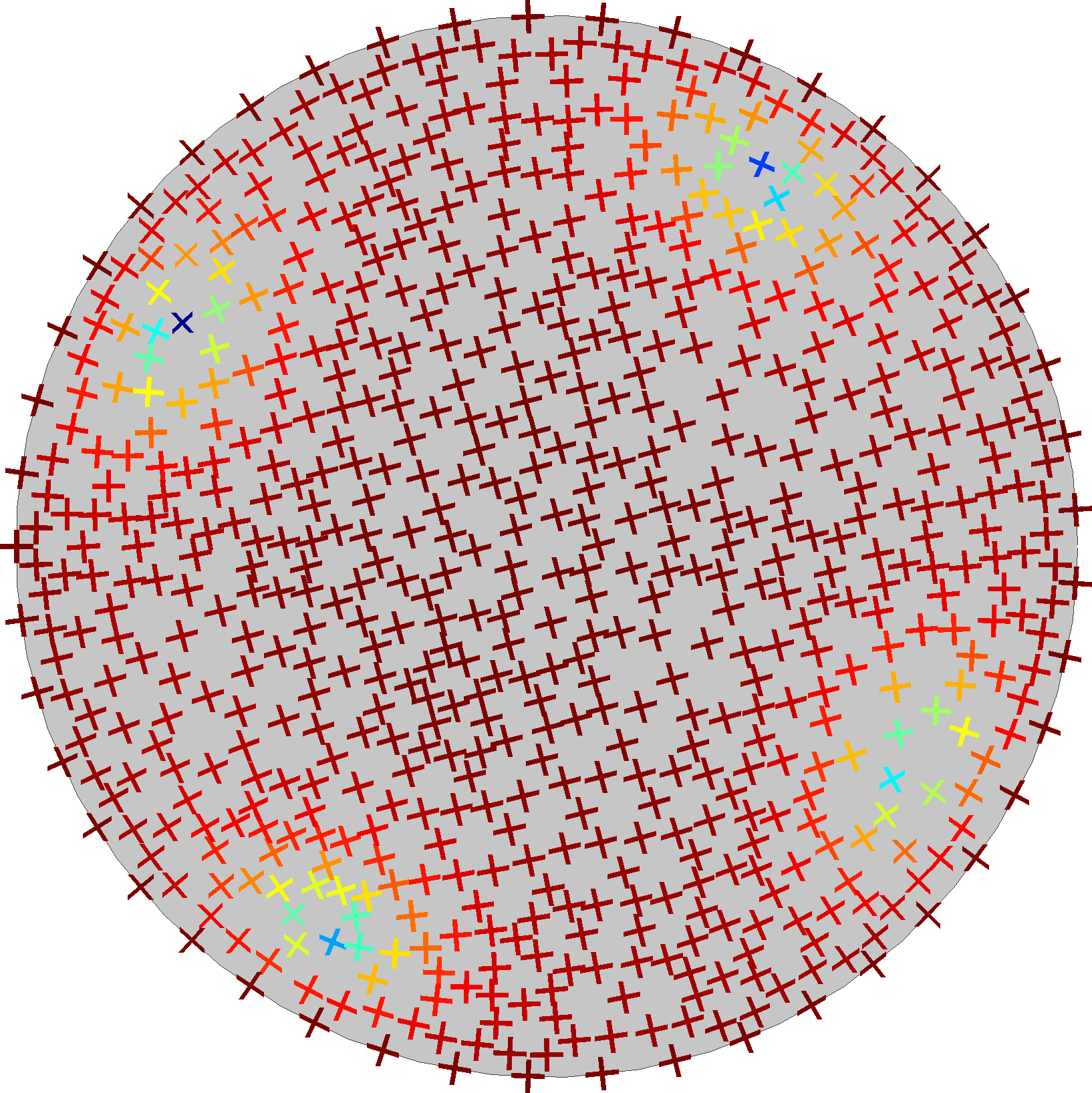}\label{sub:diskgl}}
\end{center}
\caption{Crossfield over a disk. The color describes the field norm: blue is close to zero, red close to unity.}
\label{fig:diskcf}
\end{figure}

Numerical experiments show that the norm of the vector field ${\mathbf f}$ 
computed by Laplacian smoothing (see previous section)
decreases quite rapidly as one moves away from the boundary $\partial S$,
leaving in practice large zones in the bulk of the computational domain
where the solution is small, 
and the computed crossfield inaccurate, Fig. \ref{sub:diskdir}.
A more satisfactory formulation consists of ensuring 
that the norm of ${\mathbf f}$ remains unitary
over the whole computational domain, Fig. \ref{sub:diskgl}. 
This problem can be formulated in variational form
in terms of the Ginzburg-Landau functional  
\begin{equation}
E(f_1, f_2) = 
\underbrace{{1 \over 2}\int_{\surface} \left(\left| \nabla f_1 \right|^2 + \left| \nabla
      f_2 \right|^2\right)  d\surface}_{\text{smoothing}} + 
\underbrace{{1 \over 4\epsilon^2}\int_{\surface}\left(f_1^2 + f_2^2 - 1\right)^2 d\surface}_{\text{penality}}.
\label{eq:energy2d}
\end{equation}
The first term minimizes the gradient of the crossfield 
and is therefore responsible for 
the laplacian smoothing introduced in the previous section.
The second term is a penality term that vanishes when ${\mathbf f} \in S^1$. 
The penality parameter $\epsilon$, called \emph{coherence length},
has the dimension of a length. 
The Euler-Lagrange equations of the functional \eqref{eq:energy2d}  
are the quasi-linear PDE's
\begin{equation}
\nabla^2 f_i - {1 \over \epsilon^2}\left (f_1^2 + f_2^2 - 1\right ) f_i= 0 
\quad i=1,2.
\label{eq:gl2d}
\end{equation}
called Ginzburg-Landau equations. 
If $\epsilon$ is small (enough) with respect to the dimension of $\surface$,
then ${\mathbf f}$ is of norm $1$ 
everywhere but in the vicinity of the isolated critical points ${\mathbf  x}_i^c$.

The asymptotic behavior of Ginzburg-Landau energy 
can be written as 
\begin{equation}
E = \pi \left(\sum_{i=1}^N \text{index} ({\mathbf  x}_i^c)^2 \right) \log (1/\epsilon) +
W + {\mathcal O} (1/|\log \epsilon|).
\label{eq:asygl}
\end{equation}
with
\begin{equation}
W = -\pi \sum_{i=1}^N \sum_{\underset{j\neq i}{j=1}}^N \text{index} ({\mathbf
  x}_i^c) ~\text{index} ({\mathbf  x}_j^c) \log | {\mathbf  x}_i^c - {\mathbf
  x}_j| + \mathcal{R}
\label{eq:W}
\end{equation}
as $\epsilon \rightarrow 0$  
(see \cite{bethuel1994ginzburg}, Introduction, Formulae 11 and 12).

In asymptotic regime, the energy is thus composed of three terms. 
The first term of \eqref{eq:asygl} blows up as $\epsilon \rightarrow 0$,
i.e. energy becomes unbounded if critical points are present. 
When $\epsilon$ is small, this first term dominates,
and one is essentially minimizing 
$\sum_{i=1}^N \text{index} ({\mathbf  x}_i^c)^2$
with the constraint (\ref{eq:phtimes4}). 
This indicates that a critical point of index $2$ 
has a cost of $4$ in terms of energy,
whereas 2 critical points of index $1$ have a cost of $2$. 
All critical points should therefore be of index $\pm 1$,
and their number should be $N = 4 ~ |\chi|$.
This is indeed good news for our purpose\,: 
good crossfields should have few critical points of lower indices.

The second term of \eqref{eq:asygl} 
is the \emph{renormalized energy} $W$ (\ref{eq:W}).
It remains bounded when $\epsilon$ tends to $0$. 
The double sum in $W$ reveals the existence of a logarithmic force 
between critical points.
The force is attractive between critical points with indices of opposite signs,
and repulsive between critical points with indices of the same signs. 
The second term in (\ref{eq:W}) 
is more complicated and is detailed in \cite{bethuel1994ginzburg}. 
Basically, $\mathcal{R}$ represents a
repulsing  force that forbids critical points to approach the boundaries.

Finally, the third term in \eqref{eq:asygl} vanishes as $\epsilon \rightarrow 0$. 
At the limit, all energy is thus carried by the critical points of the field. 
All this together allows to believe 
that Ginzburg-Landau model is a good choice for computing crossfields. 
It produces few critical points, which are moreover 
well-distributed over the domain. 

\section{Computation of Crossfields: Nonplanar Generalization}\label{sec:fem}

\begin{figure}[ht!]
\begin{center}
\subfloat[Global representation.]{\includegraphics{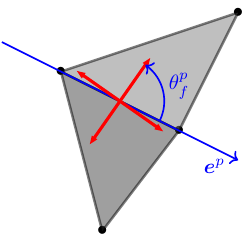}\label{sub:global}}
\subfloat[Local representation.]{\includegraphics{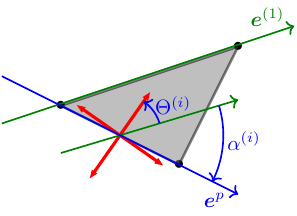}\label{sub:local}}
\end{center}
\caption{Crossfield over the $p^{\text{th}}$ edge of a mesh.}
\label{fig:cr}
\end{figure}

The finite element computation method for crossfields
is now generalized to the case of nonplanar surfaces. 
Consider the conformal triangulation
$\surface = \cup_{ijk} \Omega_{ijk}$ 
of a nonplanar surface manifold $\surface$,
each triangle $\Omega_{ijk}$ being defined 
by the vertices ${\mathbf p}_i$, ${\mathbf p}_j$ and ${\mathbf p}_k$.
Since no global reference frame exists on a nonplanar surface,
a local reference frame is associated to each edge of the triangulation.
Let $e^p$ be the $p^\mathrm{th}$ edge of the mesh,
joining nodes ${\mathbf p}_i$ and ${\mathbf  p}_j$,
and ${\mathbf n}^p$ be the average of the normals vectors
of the two triangles adjacent to $e^p$.
The vectors 
$$
\left\{ 
{\mathbf e}^p = {\mathbf p}_j-{\mathbf p}_i, 
{\mathbf t}^p = {\mathbf n}^p \times {\mathbf  e}^p
\right\}
$$
form a local frame $\{\hat{e}^p,\hat{t}^p\}$ which enables the representation of the connector values 
of the discretized crossfield ${\mathbf f}$,
$$
f^p_1 = \cos 4 \theta_f^p
\quad , \quad
f^p_2 = \sin 4 \theta_f^p,
$$
which are attached to the center of the edges of the triangulation.
Actually, $\theta_f^p$ is assumed to be the same along $\bm{e}^p$ within both planes of triangles sharing $e^p$. 
This assumption eases computations and gives a planar-like representation, Fig.~\ref{sub:global}.

\begin{figure}[ht!]
\begin{center}
\includegraphics{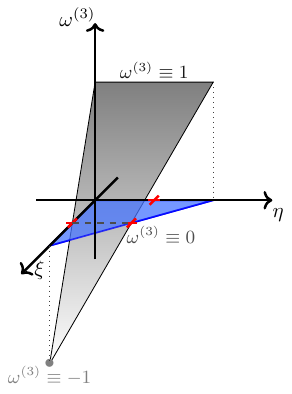} 
\end{center}
\caption{Third Crouzeix-Raviart function shape (shaded in grey) over reference triangle (in blue).}
\label{fig:crouzeix}
\end{figure}

As the connector values are attached to the edges of the mesh,
and not to the nodes, 
Crouzeix-Raviart interpolation functions are used 
instead of conventional Lagrange shape functions, \cite{crouzeix1973conforming}.
The Crouzeix-Raviart shape functions $\omega^p$
equal $1$ on corresponding edge $e^p$, and $-1$ 
on the opposite vertices (Fig.~\ref{fig:crouzeix}) in the two adjacent triangular elements.
They are polynomial
and their analytic expression in the reference triangle 
$\{\xi \in [0,1],~~ \eta \in [0,1-\xi]\}$ reads
$$
\omega^{(1)}(\xi,\eta) = 1-2\eta
\quad , \quad
\omega^{(2)}(\xi,\eta) = 2(\xi+\eta)-1
\quad , \quad
\omega^{(3)}(\xi,\eta) = 1-2\xi,
$$
where indices $(1), (2)$ and $(3)$ enclosed in parentheses
denote the local edge numbering in the considered triangular element. 
 
Each of the three edges of a triangle $\Omega_{ijk}$ has its own local reference frame.
If one is to interpolate expressions involving the vector field ${\mathbf  f}$
over this element,
the three edge-based reference frames 
have to be appropriately related with each other (see \cite{ray2016practical}).
We arbitrarily take the reference frame of the first edge of the element as reference,
and express the angular coordinate of the two other edges 
in function of this one with the relationships (Fig.~\ref{sub:local})
$$
\Theta^{(1)} = \theta_f^{(1)}
\quad , \quad
\Theta^{(2)} = \theta_f^{(2)} + \alpha^{(2)}
\quad , \quad
\Theta^{(3)} = \theta_f^{(3)} +  \alpha^{(3)}.
$$

Thus, the 6 local unknowns of triangle $\Omega_{ijk}$ 
can be expressed as a
function of the 6 edge unknowns by
{\fontsize{5}{6}\selectfont $$
\begin{pmatrix}
\cos{4 \Theta^{(1)}} \\
\cos{4 \Theta^{(2)}} \\
\cos{4 \Theta^{(3)}} \\
\sin{4 \Theta^{(1)}} \\
\sin{4 \Theta^{(2)}} \\
\sin{4 \Theta^{(3)}} \\
\end{pmatrix}
\equiv
\underbrace{
\begin{pmatrix}
F^{(1)}_1 \\
F^{(2)}_1 \\
F^{(3)}_1 \\
F^{(1)}_2 \\
F^{(2)}_2 \\
F^{(3)}_2 \\
\end{pmatrix}
}_{{\mathbf F}^{ijk}}
=
\underbrace{
\begin{bmatrix}
1 &  0                                   & 0   & 0 & 0 & 0 \\
0 &  \cos{4 \alpha^{(2)}}   & 0 & 0 & \sin{4 \alpha^{(2)}} & 0 \\
0 &  0                                   & \cos{4 \alpha^{(3)}} & 0& 0 & \sin{4 \alpha^{(3)}} \\
0 &  0                                   & 0 & 1 & 0 & 0 \\
0 &  -\sin{4 \alpha^{(2)}}  & 0 & 0 & \cos{4 \alpha^{(2)}} & 0 \\
0 &  0                                   & -\sin{4 \alpha^{(3)}} & 0
& 0 & \cos{4 \alpha^{(3)}} \\
\end{bmatrix}
}_{{\mathbf R}^{ijk}}
\underbrace{
\begin{pmatrix}
f^{(1)}_{1} \\
f^{(2)}_{1} \\
f^{(3)}_{1} \\
f^{(1)}_{2} \\
f^{(2)}_{2} \\
f^{(3)}_{2} \\
\end{pmatrix}
}_{{\mathbf f}^{ijk}}
$$}
and we have the interpolation
$$
F^{ijk}_1 (\xi, \eta) = \sum_{i=1}^3 \omega^{(i)} (\xi, \eta) \ F^{(i)}_1
\quad , \quad
F^{ijk}_2 (\xi, \eta) = \sum_{i=1}^3 \omega^{(i)}(\xi, \eta) \ F^{(i)}_2
$$
for the vector field ${\mathbf f}$ in the triangle $\Omega_{ijk}$.

A Newton scheme is proposed to converge to the solution. 
The Newton iteration at stage 
$n$ for solving \eqref{eq:gl2d} consists of solving:

\begin{equation}\label{eq:nr}
  \scriptsize
\nabla^2 \begin{pmatrix}f_{1} \\ f_{2}\end{pmatrix}_{n} 
-  
{1 \over \epsilon^2}
\begin{pmatrix} 
3 f_1^2 + f_2^2 - 1 & 2 f_1f_2 \\ 
2 f_1 f_2 & f_1^2 + 3 f_2^2 - 1
\end{pmatrix}_{n-1}
\begin{pmatrix}f_{1} \\ f_{2}\end{pmatrix}_{n}
=
- {2 \over \epsilon^2}
\begin{pmatrix} 
f_1^3 + f_1 f_2^2 \\
f_1^2 f_2 + f_2^3
\end{pmatrix}_{n-1}
\end{equation}
The $6 \times 6$ elementary matrix ${\mathbf K}^{ijk} $ 
and the $6\times 1$ elementary vector ${\mathbf B}^{ijk}$ 
of element $\Omega_{ijk}$ 
are then given by ${\mathbf K}^{ijk}$ 
{\fontsize{4}{8}\selectfont
  \begin{equation}  
\begin{bmatrix}
\left(\bigints_{\Omega_{ikj}} \nabla \omega_m \cdot \nabla \omega_n +
  {1 \over \epsilon^2} (3F_1^2 + F_2^2 - 1) \omega_m \omega_n ~ d\Omega \right)
& 
\left(\bigints_{\Omega_{ijk}} {2 \over \epsilon^2} F_1 F_2 \omega_m \omega_n ~ d\Omega \right)\\
\left(\bigints_{\Omega_{ijk}} {2 \over \epsilon^2} F_1 F_2 \omega_m \omega_n ~ d\Omega \right)&
\left(\bigints_{\Omega_{ikj}} \nabla \omega_m \cdot \nabla \omega_n +
   {1 \over \epsilon^2} (F_1^2 + 3F_2^2 - 1)  \omega_m \omega_n ~ d\Omega \right)
\end{bmatrix}
\label{eq:nrfem}
\end{equation}}
and
\begin{equation}
{\mathbf B}^{ijk}  =  \begin{pmatrix}
\left(\bigints_{\Omega_{ijk}} \nabla F_1 \cdot \nabla \omega_n + {1 \over \epsilon^2} (F_1^3+F_2^2-F_1) \omega_n ~ d\Omega \right)\\
\left(\bigints_{\Omega_{ijk}} \nabla F_2 \cdot \nabla \omega_n + {1 \over \epsilon^2} (F_1^2+F_2^3-F_2) \omega_n ~ d\Omega \right)
\end{pmatrix}.
\label{eq:nr2}
\end{equation}
with $m,n=1 \dots 3$.

It is then necessary to transform those elementary matrix and vector
in the reference frames of the edges as
$$ {\mathbf k}^{ijk}   =   ({\mathbf R}^{ijk})^T {\mathbf K}^{ijk}
{\mathbf R}^{ijk} ~~~\text{and}~~~{\mathbf b}^{ijk}= ({\mathbf
  R}^{ijk})^T {\mathbf B}^{ijk}. $$
Then, standard finite element assembly can be performed. 
Boundary conditions are simply 
$$
f^p_1=1
\quad , \quad
f^p_2=0
$$ 
on every edge $e^p$ of $\partial \surface$. 
This nice simplification is due to the fact that unknowns are defined on
the reference frame of the edges.

\section{Numerical Validation: the Unit Disk}\label{sec:disk}

We compute the analytical location of critical points of a directional field defined on the unit disk.
The calculations are based on the Ginzburg-Landau results, described in section \ref{sec:glm}.
The numerical location obtained by our FEM is compared to the analytical one.

Let $\surface$ be the open unit disk in $\Re^2$, i.e. 
$$
\surface := \left\{(x_1;x_2)\in\Re^2 ~|~ x_1^2 + x_2^2 < 1 \right\}
$$

For a star-shaped planar domain such as $\surface$ with a smooth boundary $\partial \surface$ of exterior normal $\nu$ and tangent $\tau$,  whose vector field has $d$ critical points of index $+1$ at $X^{c} = \{\mathbf{x}_1^c,...,\mathbf{x}_d^c\}\in  \surface$, the asymptotic energy $E_\epsilon$ (in complex form) becomes
\begin{equation}\label{eqn:asymptoticGL}
E_\epsilon \underset{\epsilon \rightarrow 0}{\longrightarrow} \pi d ~|\log(\epsilon)| + W(X^c) + \mathcal{O}(\epsilon)
\end{equation}
where $W(X^c)$ is the renormalized energy
\begin{equation}\label{eqn:renormalized}
W(X^c) = -\pi \sum_{i\neq j} \log|\mathbf{x}_i^c-\mathbf{x}_j^c| + {1 \over 2} \int_{\partial \surface} \Phi ~ f \times \nabla f \cdot \mathbf{\tau} ~ ds - \pi \sum_i R(\mathbf{x}_i^c)
\end{equation}
where $\Phi$ is given by the following Neumann problem
\begin{equation}\label{eq:phineumann}
\left.\begin{array}{rclcc}
\nabla^2 \Phi(\mathbf{x}) &=& \displaystyle{2\pi \sum_{i=1}^d \delta(\mathbf{x}-\mathbf{x}_i^c)} &\text{in} & \surface \\ \\
\nabla \Phi \cdot \mathbf{\nu} &=& f \times \nabla f \cdot \mathbf{\tau} & \text{on} & \partial \surface
\end{array}\right\}
\end{equation}
and $R$ is the regular part of $\Phi$:
\begin{equation}
R(\mathbf{x}) = \Phi(\mathbf{x}) - \sum_{i=1}^d \log|\mathbf{x}-\mathbf{x}_i^c|
\end{equation}

$E_\epsilon$ is minimum when the critical points are located appropriately, i.e. when (\ref{eqn:renormalized}) is minimum. 
The renormalized energy $W$ corresponds to the Ginzburg-Landau energy (\ref{eqn:asymptoticGL}) when the singular core energy $\pi d |\log(\epsilon)|$ has been removed.
Since $W$ depends only on the location of the critical points, it is possible to compute their location in the case of the unit disk, in order to get an optimal directional field.

The minimum of $W$ is obtained by sampling points within the unit disk.
It is assumed that the $d$ critical points exhibit the $d$ symmetries of their group (the quadrilateral group in the case $d=4$).
In other words, it means that they are at the same distance $r^c$ from the center of the disk (i.e. the origin $(0;0)$), and separated two-by-two with an angle of $2\pi/d$ radians.

The Neumann problem (\ref{eq:phineumann}) is solved  by decomposing $\Phi = \Phi^0 + \Phi^1$.
The first term $\Phi^0$ is the Green function of a two-dimensional Laplacian operator, while the second one $\Phi^1$ is obtained by separation of variables $(r,\theta)$.
The solution is then
\begin{equation}
\Phi(r,\theta) =  \sum_{i=1}^d \left[  \underbrace{\vphantom{\sum_{n=1}^\infty} \frac{1}{2} \log (r^2-2r~r^c~\cos(\theta) + {r^c}^2)}_{\Phi_i^0} + \underbrace{\sum_{n=1}^\infty A_{i,n} r^n \cos(n~\theta)}_{\Phi_i^1} \right]
\end{equation}
where $A_{i,n}$ depends on the location of the i-th critical point, which is parameterized by $r^c$.
It is possible to show that the second term of (\ref{eqn:renormalized}) is zero, Appendix \ref{app:integral}. The analytical solution of Neumann problem is derived into the Appendix \ref{app:analytic}.

\begin{figure}
\begin{center}
\subfloat[Four critical points: $r^{c^*}=0.85$.]{\includegraphics[width=.45\linewidth]{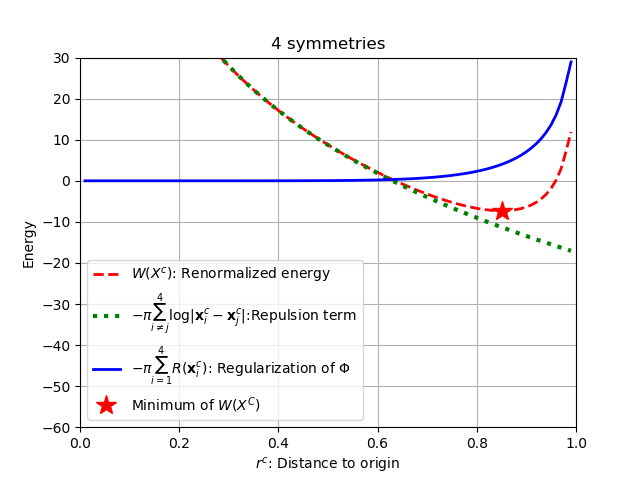}}\hspace*{.25cm}
\subfloat[Six critical points: $r^{c^*}=0.90$.]{\includegraphics[width=.45\linewidth]{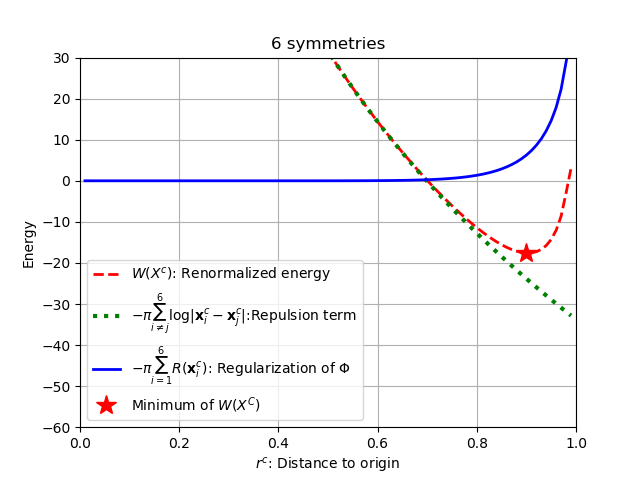}}
\caption{Python evaluations of renormalized energy $W$ for different $X^c$ on a unit disk $\surface$.}\label{fig:aloce}
\end{center}
\end{figure}

The evaluation of $W$ consists of computing the first and last terms, by sampling the disk.
The sampling is done by selecting $d$ critical points spaced by $2\pi/d$ radians.
The distance $r^c$ is sampled between zero and one.
The distance $r^{c^*}$ which gives the lowest value of $W$ defines the location of the critical points.
A Python script performs the evaluations and returns the optimal distance $r^{c^*}$, Fig. \ref{fig:aloce}.

\begin{figure}
\begin{center}
\begin{tikzpicture}
\node at (-.25\linewidth,0) {
\subfloat[Crossfield.]{\includegraphics[width=.45\linewidth]{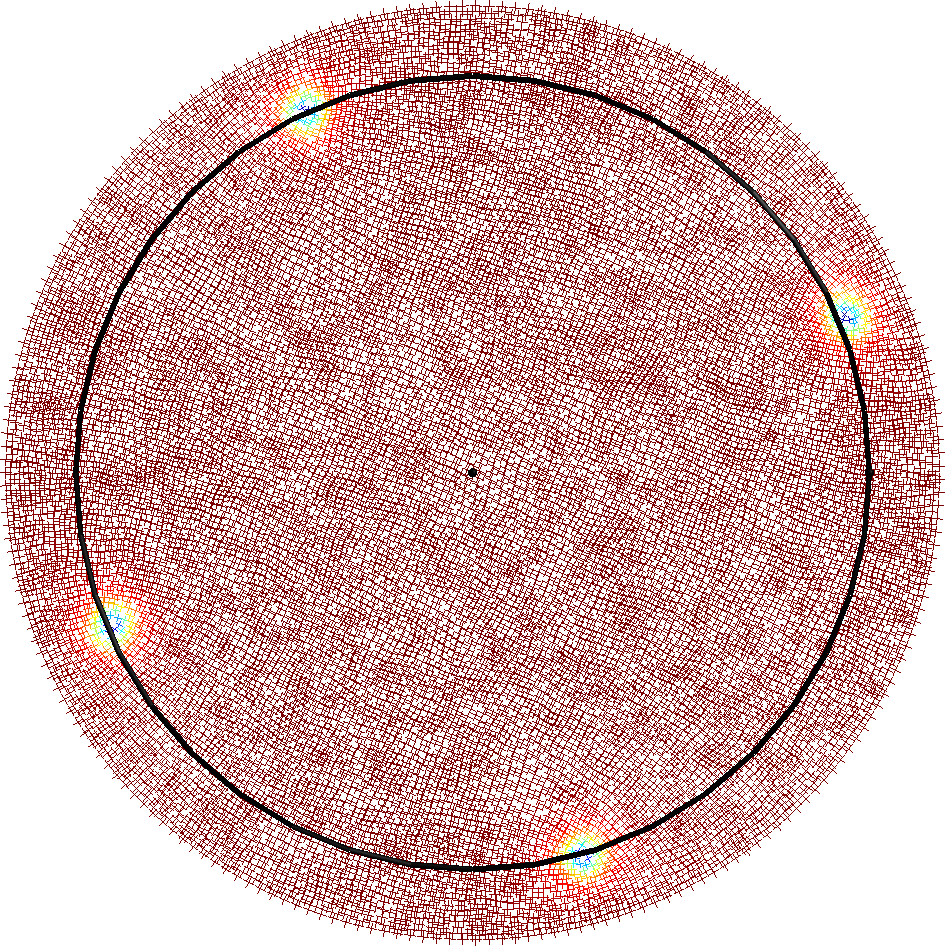}\label{sub:adcf}}
};
\node at (.25\linewidth,0) {
\subfloat[Asterisk field.]{\includegraphics[width=.45\linewidth]{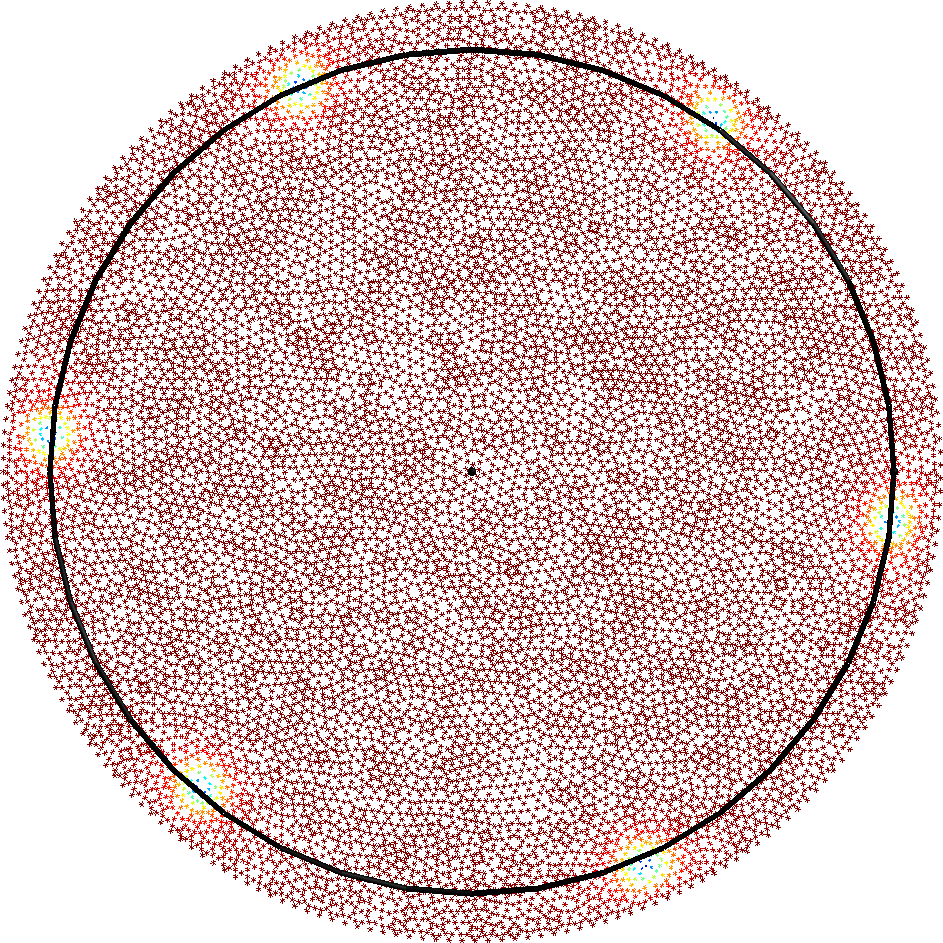}\label{sub:adaf}}
};
\draw[very thick,<->] (-.25\linewidth,.225)--(-.0587\linewidth,.225);
\node at (-.15\linewidth,.4) {$0.85$};
\draw[very thick,<->] (.25\linewidth,.05)--(.0475\linewidth,.05);
\node at (.15\linewidth,.3) {$0.90$};
\end{tikzpicture}
\caption{FEM computations of direction fields on a unit disk $\surface$: the critical points are in blue areas.}\label{fig:aloc}
\end{center}
\end{figure}

The corresponding directional fields are computed, and their critical point locations are compared with circles which radii correspond to $r^{c^*}$, Fig. \ref{fig:aloc}.
The location of critical points are really close to the estimation based on the analytical solution of $W(X^c)$ in the case of the unit circle.
They tend to draw the corners of the polygon of symmetry: a square in the case of the crossfield, Fig. \ref{sub:adcf} and a regular hexagon for the asterisk field, Fig. \ref{sub:adaf}.
The critical points are quite close to the unit circle. 
The more critical points, the closer to the unit circle they are. 
We understand that the repulsion term is stronger than the regularization term within the domain.
The regularization term is only able to forbid critical points to be on the boundary, i.e. the unit circle.

\section{A Surprising Result: the Sphere}\label{sec:sphere}

\begin{figure}[!ht]
\begin{center}
\includegraphics[width=.45\linewidth] {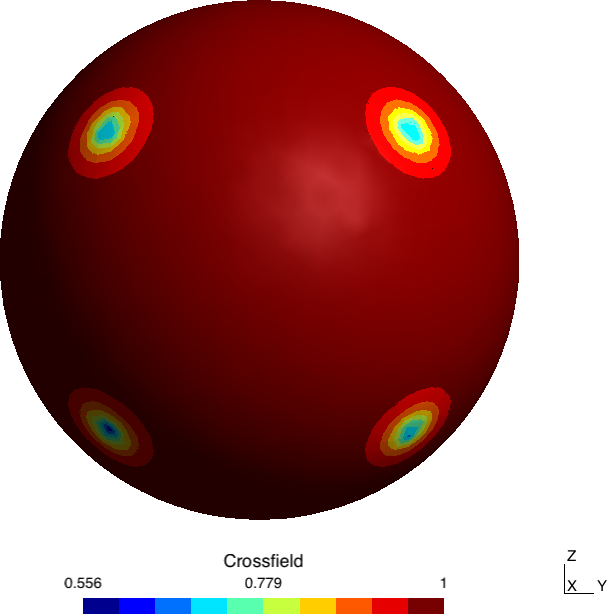}
\includegraphics[width=.45\linewidth] {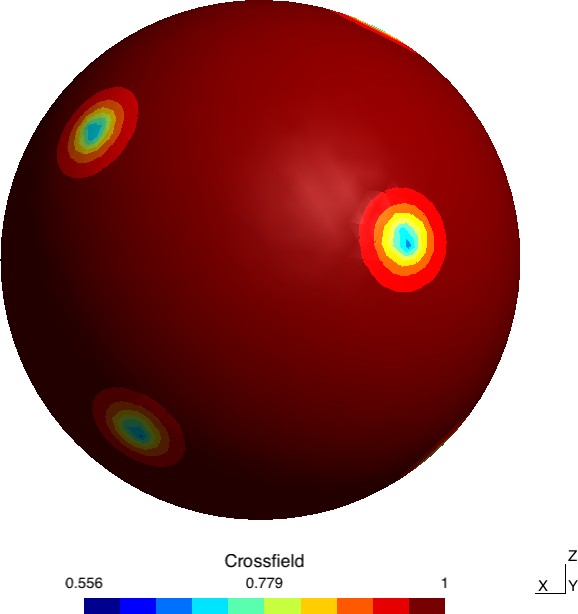}
\includegraphics[width=.45\linewidth] {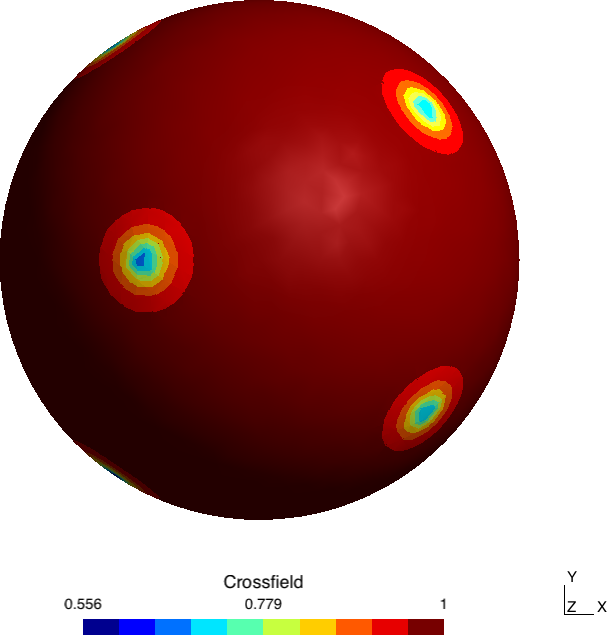}~\\
\end{center}
\caption{Mesh of the sphere. Colors correspond to the 2-norm
  $\|{\mathbf f}\|$ of the crossfield. The 8 critical points are 
  located on two squares of side $1/\sqrt{3}$, which corresponds to the
  size of the inscribed cube. The two squares are tilted by 45
  degrees. \label{fig:sphere1}}
\end{figure}

\begin{figure}[!ht]
\begin{center}
\includegraphics[width=.5\linewidth]{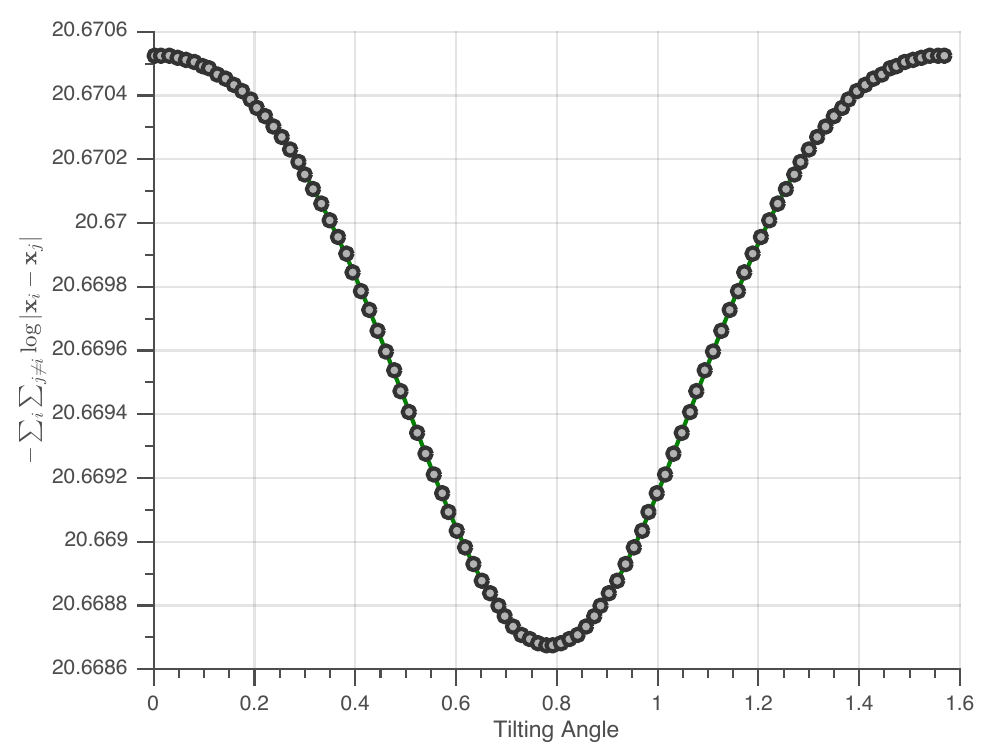}
\end{center}
\caption{Energy vs. tilting angle for the sphere. The minimum
  corresponds to a tilting angle of $\pi/4$. \label{fig:sphere2}}
\end{figure}

\begin{figure}[!ht]
\begin{center}
\includegraphics[width=.45\linewidth] {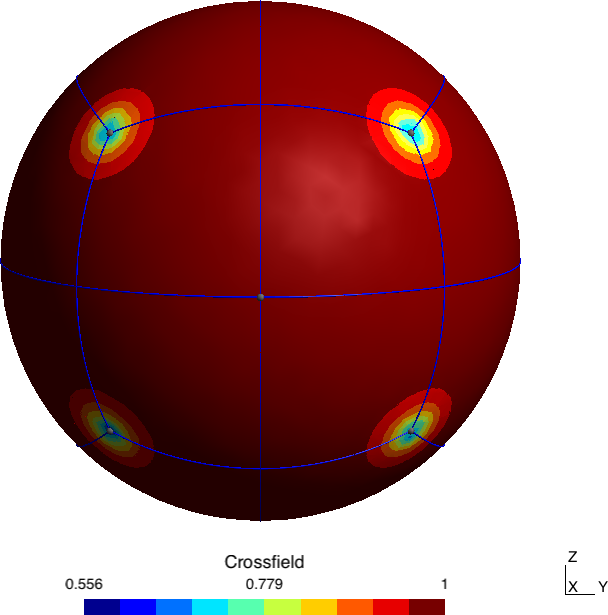}
\includegraphics[width=.45\linewidth] {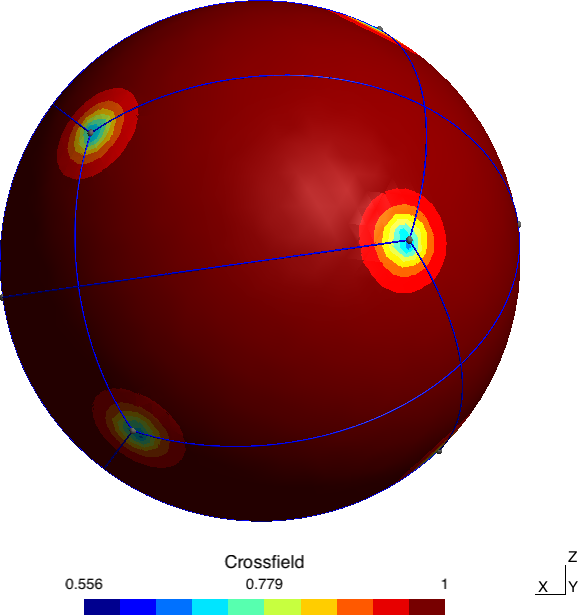}
\includegraphics[width=.45\linewidth] {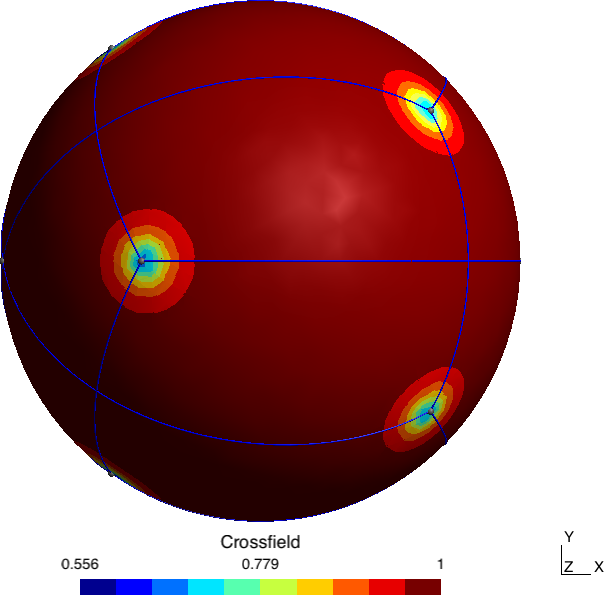}
\end{center}
\caption{Separatrices from crossfield. \label{fig:sphere3}}
\end{figure}

Let us compute the crossfield on a unit
sphere. The sphere has no boundary so we choose randomly one edge of
the mesh and fix the crossfield for this specific edge. The mesh of
the sphere is made of 2960 triangles (see Fig. \ref{fig:sphere1}).
A value of $\epsilon=0.1$ was chosen for the computation. A total of $29$ Newton iterations
were necessary to converge, by reducing the residual norm to $10^{-12}$.
The location of the $8$ critical points is indeed not what we expected: our
initial intuition was that critical points would be located at the
corners of an inscribed cube of side $1/\sqrt{3}$. In all our
computations i.e. while changing the mesh and $\epsilon$, critical
points are located on two squares of side $1/\sqrt{3}$, those two
squares being tilded by $45$ degrees around their common axe (see
Fig. \ref{fig:sphere1}). Equilateral triangle patterns are formed
between critical points that belong to both squares. 
In reality, our solution is the right solution. In the asymptotic
regime, the location ${\mathbf x}_i^c$ of the $8$ critical points tends
to minimize $-\sum_{i} \sum_{i\neq j} \log |{\mathbf x}_i^c - {\mathbf
  x}_j^c|$ (see Equations \eqref{eq:asygl} and \eqref{eq:W}). We have
thus computed $-\sum_{i} \sum_{i\neq j} \log |{\mathbf x}_i^c - {\mathbf
  x}_j^c|$ for tilting angles ranging from $0$ to $\pi/2$. Fig.
\ref{fig:sphere2} shows clearly that the minimum of the energy
corresponds to an angle of $\pi/4$, which is exactly what is found by
the finite element formulation.
Fig. \ref{fig:sphere3} shows the crossfield as well as
the separatrices. The separatrices were computed ``by hand''.
The solution that has been found is related to what is called
the Whyte's problem (cf. \cite{saff1997distributing,dragnev2002separation}) that consists of finding $N$ points
on the sphere which positions maximize the product of their
distances. The critical points are called \emph{logarithmic extreme points}
or {\textit elliptic Fekete points} (see \cite{fekete1923verteilung}).

The specific configuration that corresponds to $N=8$
is called an anticube (or square antiprism) and is exactly the one
that was found numerically.

In the case of an asterisk field, the critical points are the summits of an icosahedron, 
which is the solution of Whyte's problem for $N=12$.
This superb result shows that it is indeed possible to
use crossfields not only for building quadrangles but also to build
equilateral triangles. 

\begin{figure}[!ht]
\begin{center}
\includegraphics[width=.45\linewidth] {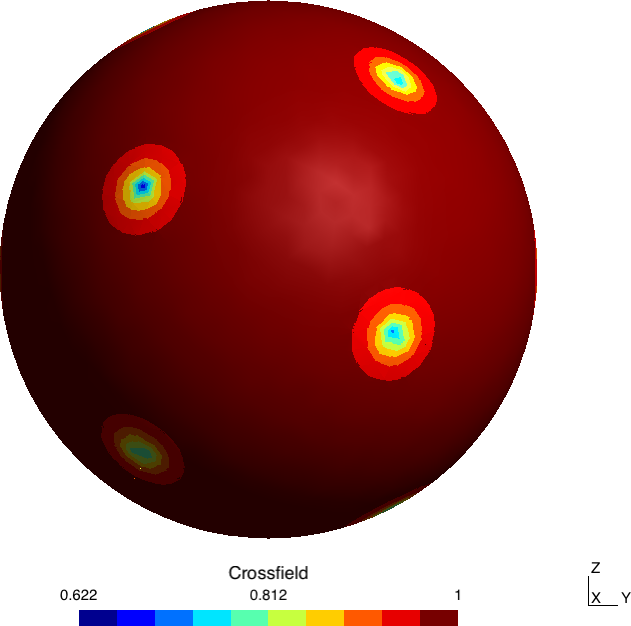}
\includegraphics[width=.45\linewidth] {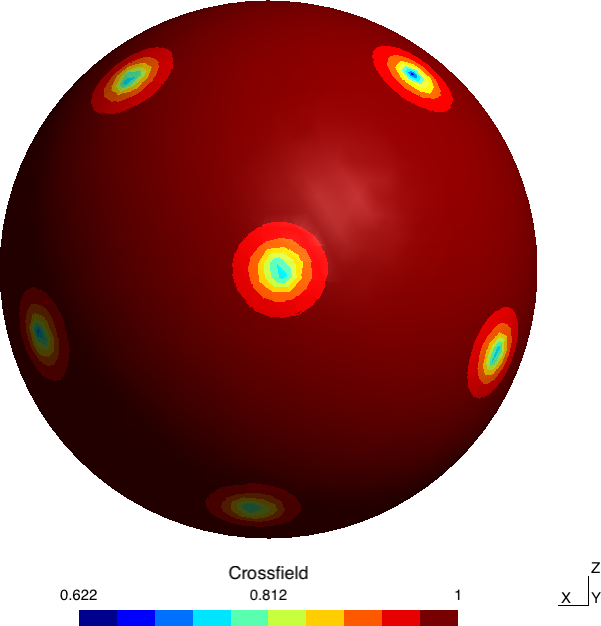}
\includegraphics[width=.45\linewidth] {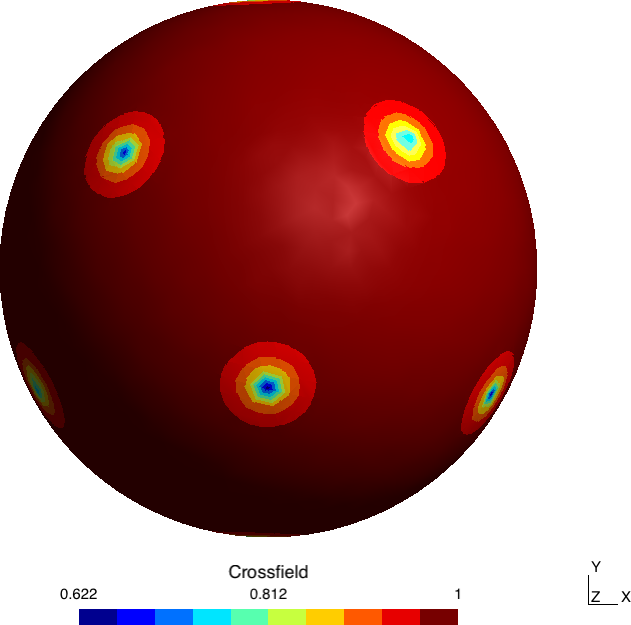}
\end{center}
\caption{Asterisk field (6 symmetries) which the critical points correspond to the
  corners of an icosahedron. \label{fig:sphere4}}
\end{figure}

Actually, it is possible to show that the critical points computed over the sphere by Ginzburg-Landau correspond to the solution of Whyte's problem for any even value of $N$ (see \cite{jezdimirovic2017elliptic}).

\section{Weak Boundary Conditions}\label{sec:naca}
\begin{figure}[!ht]
\begin{center}
\includegraphics[width=.45\linewidth]{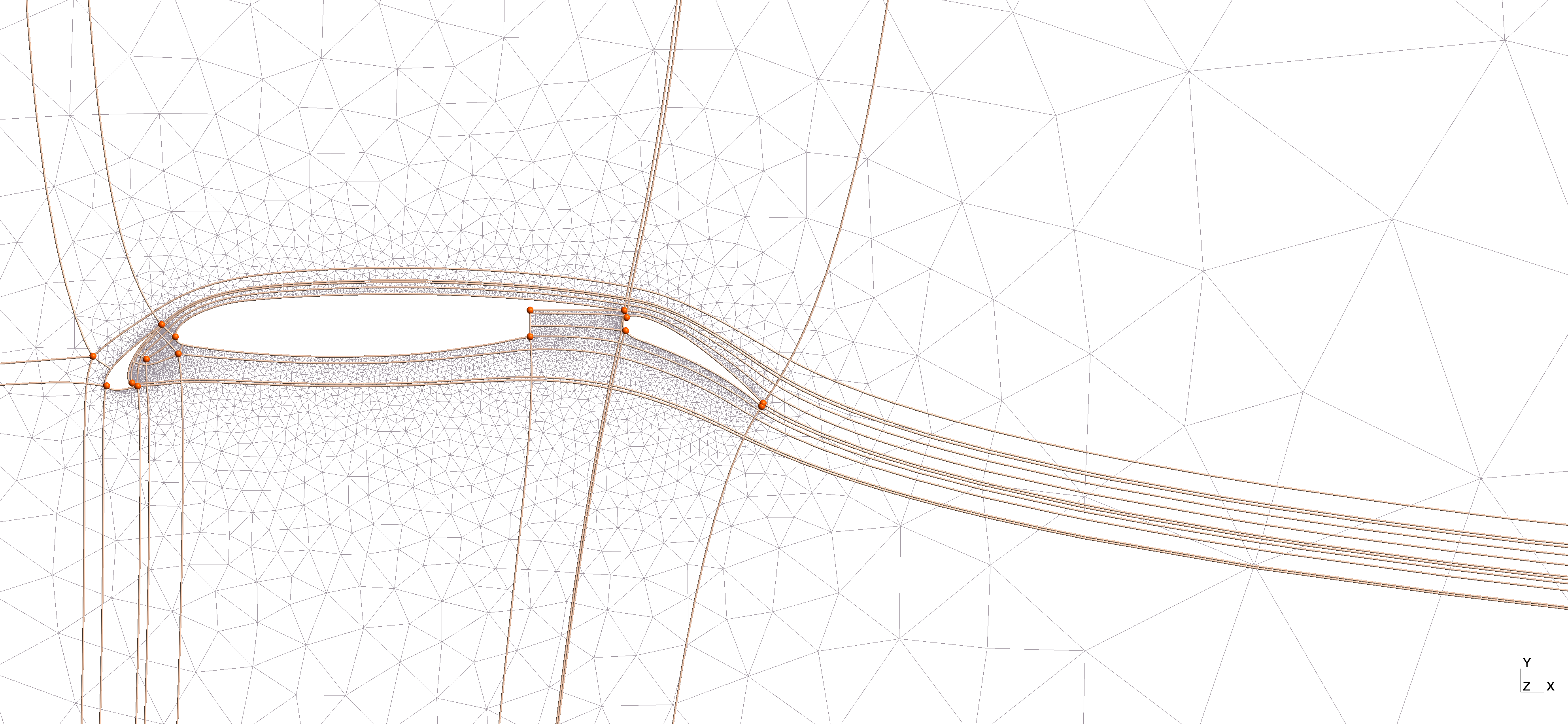}\hspace*{.5cm}
\includegraphics[width=.45\linewidth]{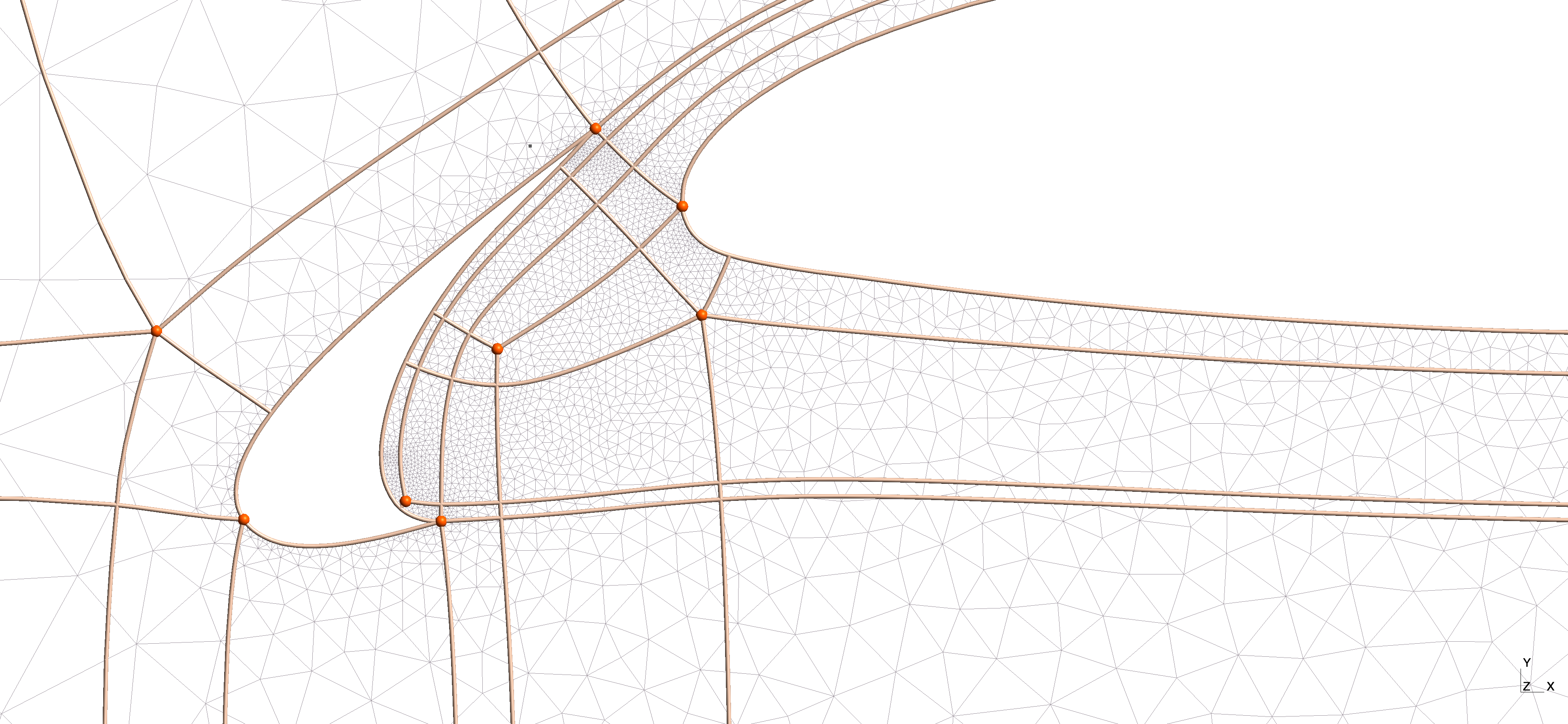}
\end{center}
\caption{Graph of singularities for the three component wing. Right
  figure is a zoom on the leading edge slat. \label{fig:wing}}
\end{figure}

In this section, we have computed the graph of singularities of a
standard CFD test case: a three component wing domain with $\chi =
-2$. This example is very similar to the one presented by \cite {kowalski2013pde}. The solution has been computed on a
non uniform triangular mesh of about $15,000$ triangles. The graph of
singularities has been depicted on Fig. \ref{fig:wing}.
Weak boundary conditions have been applied to the different components
of the wing where a penalization replaces the strong imposition of $f$
on boundaries. A new term is thus added to Energy \eqref{eq:energy2d}
for taking into account boundary conditions:
\begin{equation}
\tiny
  {1 \over 2}\int_{\surface} \left(\left| \nabla f_1 \right|^2 + \left| \nabla
      f_2 \right|^2\right)  d\surface + 
{1 \over 4\epsilon^2}\int_{\surface}\left(f_1^2 + f_2^2 -
    1\right)^2 d\surface
+ {L \over 2\epsilon^2}\int_{\partial \surface} \left[(f_1-\bar{f}_1)^2 +
  (f_2-\bar{f}_2)^2\right] d\partial \surface
\label{eq:energy2dweakbc}
\end{equation}
where $\bar{f}_1$ and $\bar{f}_2$ are values of the crosses that are
weakly imposed on the boundary and $L$ the characteristic size of the
problem. 
This new treatment allows singularities to migrate on the boundary, making their repulsive
action finite. Figure \ref{fig:wing} clearly shows that effect: a
singularity of index $1/4$ sits on the leading edge of the slat, allowing a
 clean decomposition of the domain. The same migration is also
observed on the leading edge of the profile.  A strong imposition of
boundary conditions naturally leads to singularities that are very
close to regions of the boundary with high curvature, usually at a
distance from the boundary that is one mesh size. Artificial 
boundary layers are thus added to the decomposition (see \cite[Fig. 12 and 14]{kowalski2013pde}). 


\section{Application of our FEM Scheme to the Torus}\label{sec:torus}

\begin{figure}
\begin{center}
\subfloat[Clipped view.]{\includegraphics[width=.3\linewidth]{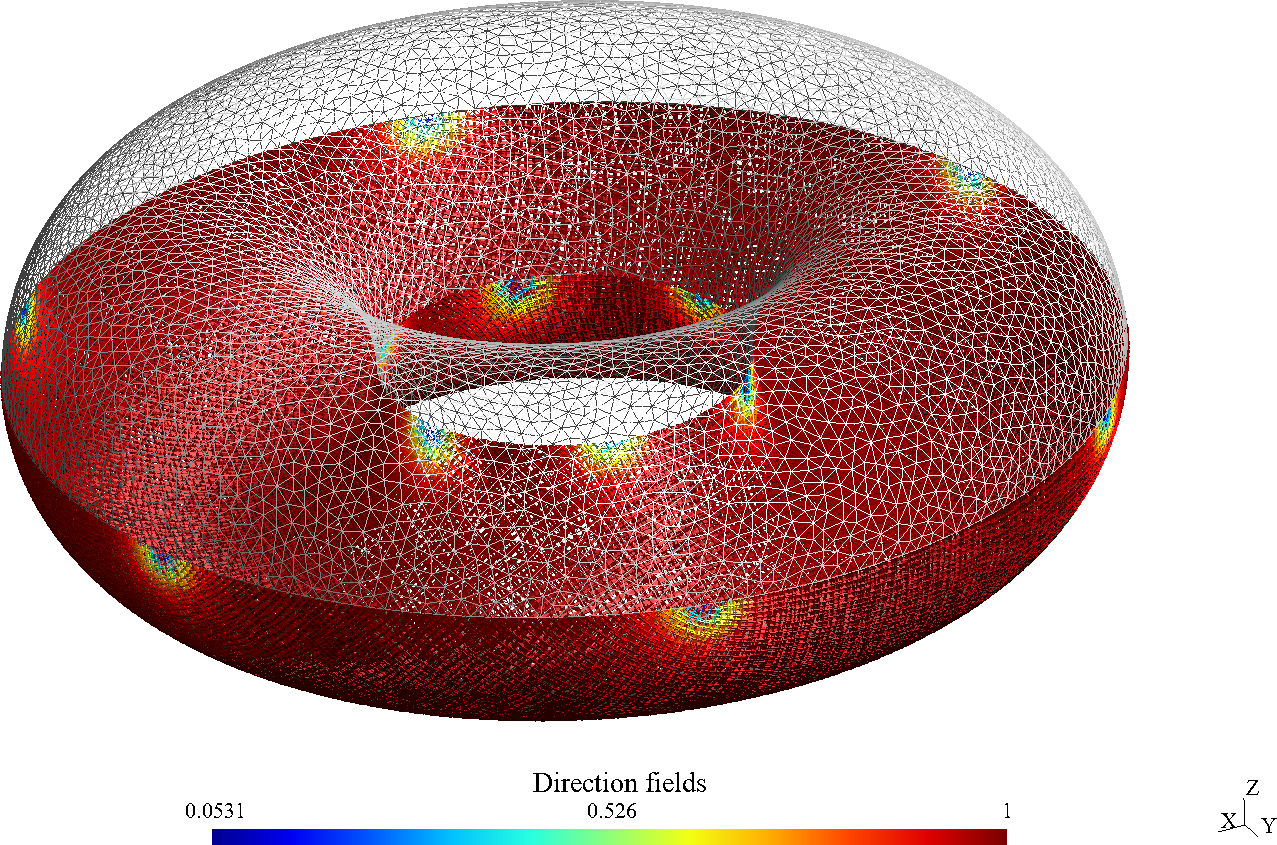}\label{sub:tcritical}}
\subfloat[Crossfield.]{\includegraphics[width=.3\linewidth]{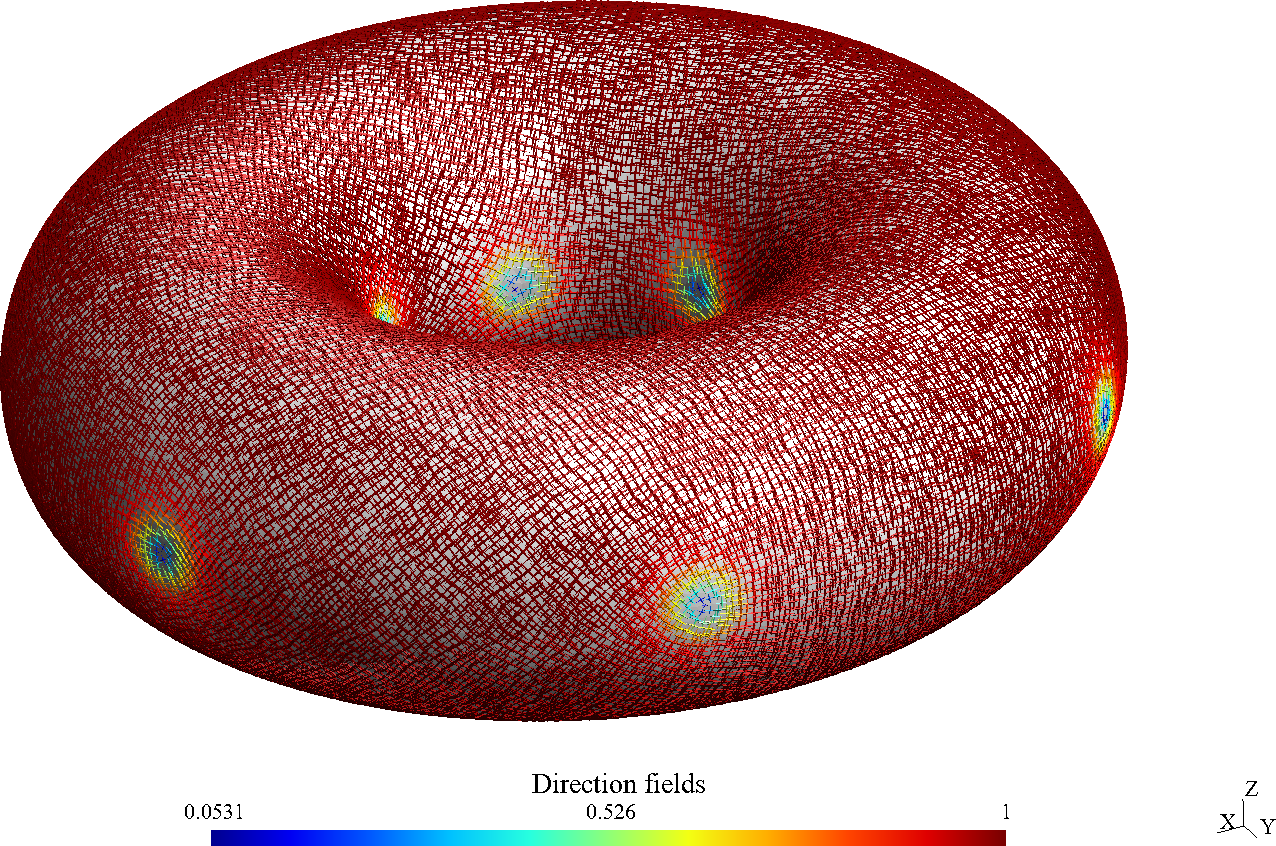}\label{sub:tcf}}
\subfloat[Polyquad decomposition.]{\includegraphics[width=.3\linewidth]{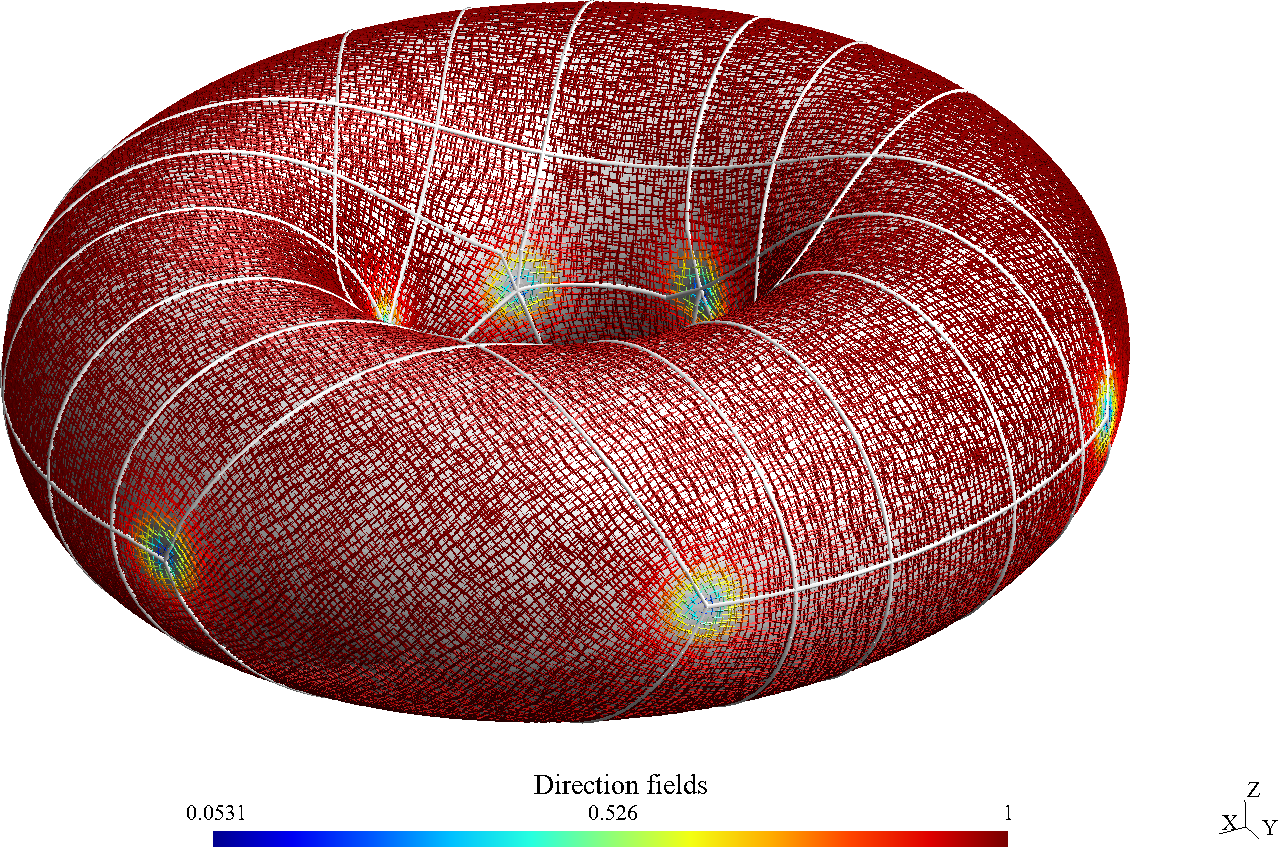}\label{sub:tquad}}
\end{center}
\caption{FEM ($\epsilon = 0.06$) crossfield on a torus discritized by 20612 triangles. Twelve critical points of opposite indices ($\pm 1/4$).}\label{fig:torus_cross}
\end{figure}

The Euler-Poincaré characteristic of the torus is $\chi = 0$. 
Theoretically, we should obtain a crossfield without critical points.
But our FEM scheme gives crossfield with twelve critical points,
located where the Gaussian curvature is maximal (exterior) or minimal (interior), 
Fig.~\ref{fig:torus_cross}.
Fig.~\ref{sub:tcritical} shows that the six critical points 
located on the maximal Gaussian curvature line  
are facing the six corresponding critical points
located on the minimal Gaussian curvature line. 
Moreover, as the former have an index $+1/4$, and the latter an index $-1/4$, Fig.~\ref{sub:tcf}, 
the index sum of the surface is zero, as predicted by the Poincaré-Hopf theorem.

Our FEM scheme does not reach however  
the asymptotic behavior $(\epsilon \rightarrow 0)$ of the Ginzburg-Landau functional.
It means that our penalty factor $\epsilon$ is not low enough.
Otherwise, the computed crossfield should not have any critical points owing to (\ref{eq:asygl}).
Actually, the computed crossfield has a lower energy ($72.10$) 
than the crossfield with no critical point that could be drawn by aligning crosses 
with the main curvatures of the surface ($84.58$).
The tentative polyquad decomposition shown in Fig.~\ref{sub:tquad}
indicates that the field computed with the Ginzburg-Landau approach tends 
to be more uniform, in order to reduce the Dirichlet energy.
It confirms that the Dirichlet term is stronger than the penalty term.

\section{Conclusion}

This article has demonstrated the consistency of the Ginzburg-Landau theory 
to compute directional fields on arbitrary surfaces. 
The proposed approach relies on a physical and mathematical backgrounds.
This provides proofs, analytical solutions
and helps delineating fundamental mathematical properties that can be exploited in algorithms. 

In particular, the Ginzburg-Landau theory states that 
when the coherence length $\epsilon$ is small enough, 
the asymptotic behavior is reached,
i.e., the number of critical points of the crossfield is minimal, 
their index is also minimal and they are optimally distributed.
A simple FEM scheme has been implemented to validate numerically this assertion.
Crossfields have been computed on the unit disk
and solutions conform with the Ginzburg-Landau theory have been found. 
The location of computed critical points on the 2-sphere corresponds to the solution of Whyte's problem: 
for a crossfield they are at the summits of an anticube 
whereas for an asterisk field they are at the summits of a regular dodecahedron.

By weakening the boundary conditions of the Ginzburg-Landau problem,
critical points are  no longer repelled in the interior of the domain and can be located on the boundary, 
which improves the polyquad decomposition in the case of the NACA profiles.



\begin{figure}
\begin{center}
\begin{tikzpicture}
\node at (0,0) {\includegraphics[width=.9\linewidth] {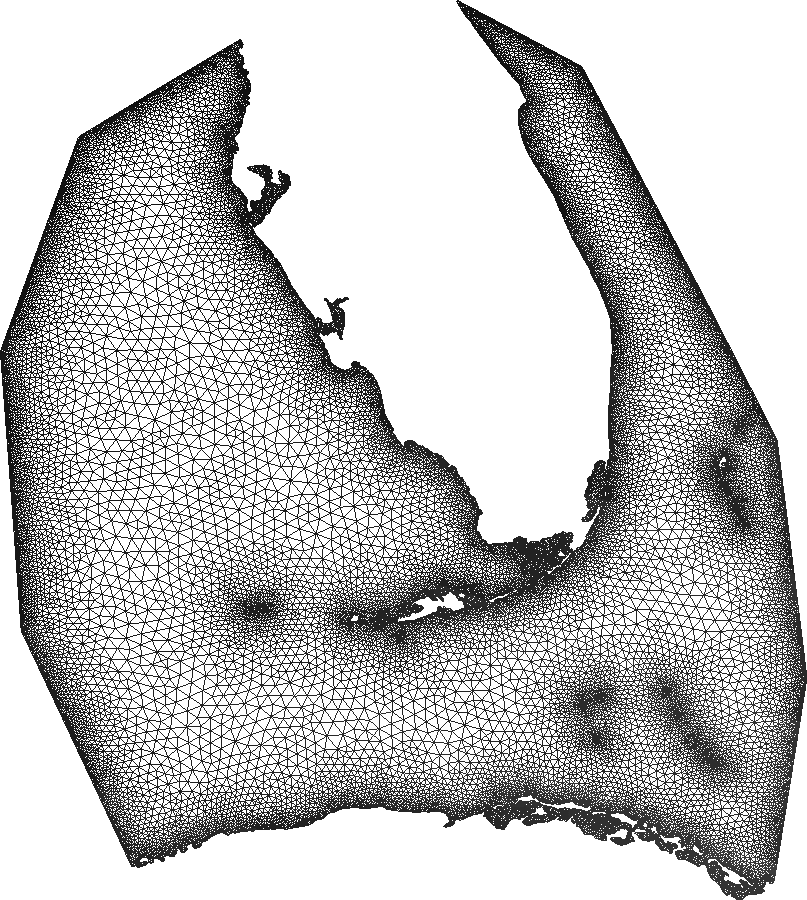}};
\draw[black!75!white,dashed,very thick] (2.25,-1.) rectangle (3.75,-.125);
\end{tikzpicture}
\end{center}
\caption{Florida keys: input triangular mesh ($\bar{\gamma}=0.9504$). The rectangle will be enlarged.} \label{fig:floridaIn}
\end{figure}

\begin{figure}[!ht]
\begin{center}
\subfloat[Right-angled triangles from crossfield.]{\includegraphics[width=.9\linewidth]{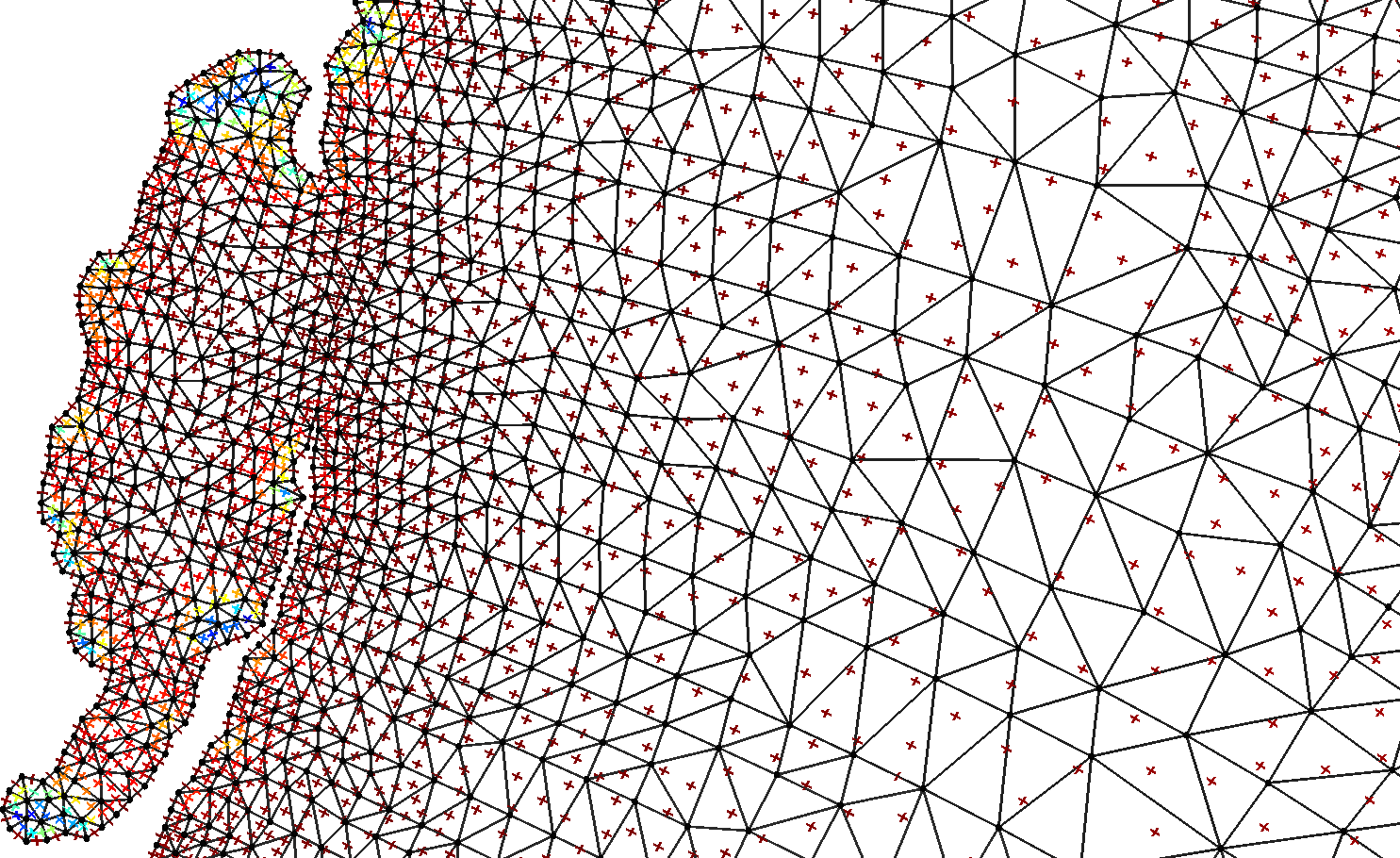}\label{sub:right}}
\hfill
\subfloat[Final quadrangular mesh after optimization.]{\includegraphics[width=.9\linewidth] {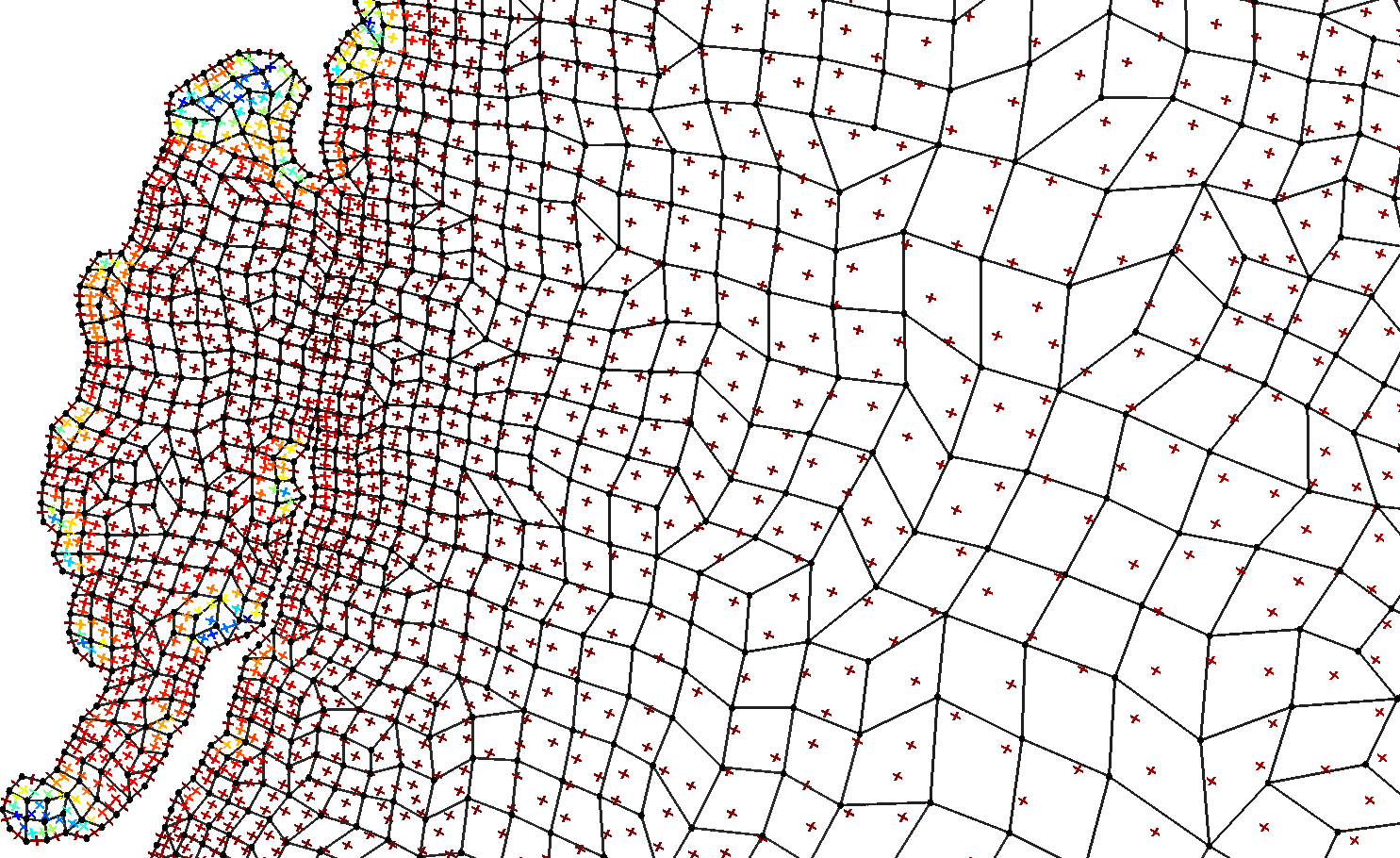}\label{sub:quad}}
\end{center}
\caption{Zoom on the Florida keys, the color map is 0 (blue) to 1 (red) and describes the norm of directions. \label{fig:floridaQ}}
\end{figure}

\begin{figure}[!ht]
\begin{center}
\includegraphics[width=.9\linewidth]{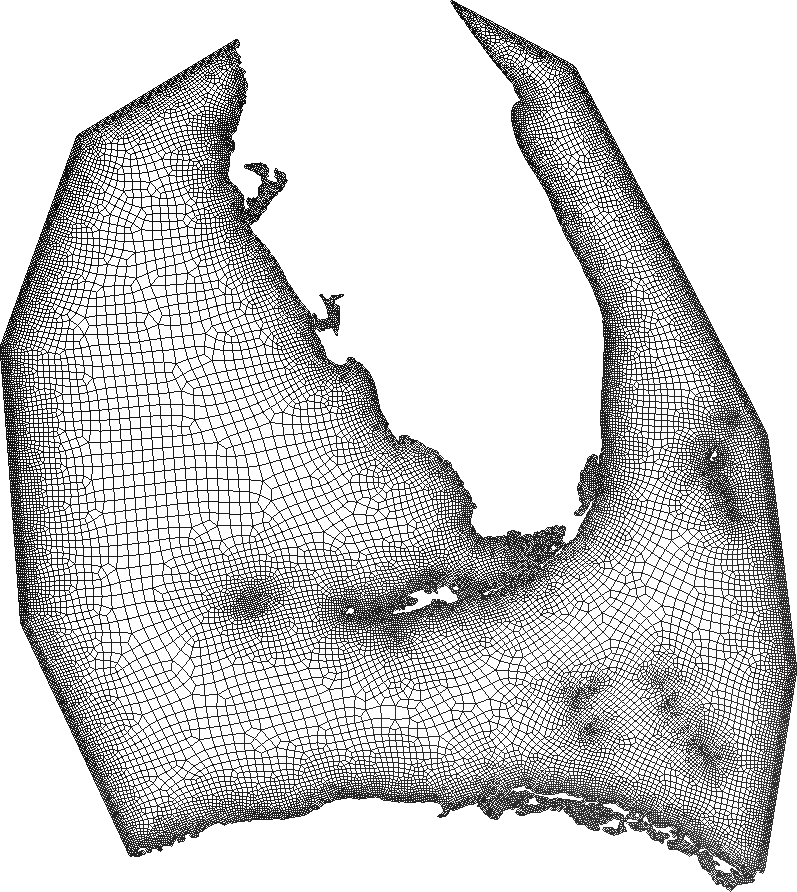}
\end{center}
\caption{Final quadrangular mesh over the Florida keys.}\label{fig:floridaOutQuad}
\end{figure}

Finally, the process is applied to the quadrangular meshing 
of the coastal domain around Florida peninsula, Fig.~\ref{fig:floridaIn}. 
Quadrangles are merged from right-angled triangles whose vertices have been spawned 
along the integral lines of  a crossfield, Fig.~\ref{sub:right}. 
One sees on Fig.~\ref{sub:quad} how the edges of the recombined quadrangular elements 
tend to follow the crossfield,  
and the final mesh is of satisfying quality, Fig.~\ref{fig:floridaOutQuad}. 

\begin{figure}[!ht]
\begin{center}
\subfloat[Asterisk field over the new triangles.]{\includegraphics[width=.9\linewidth] {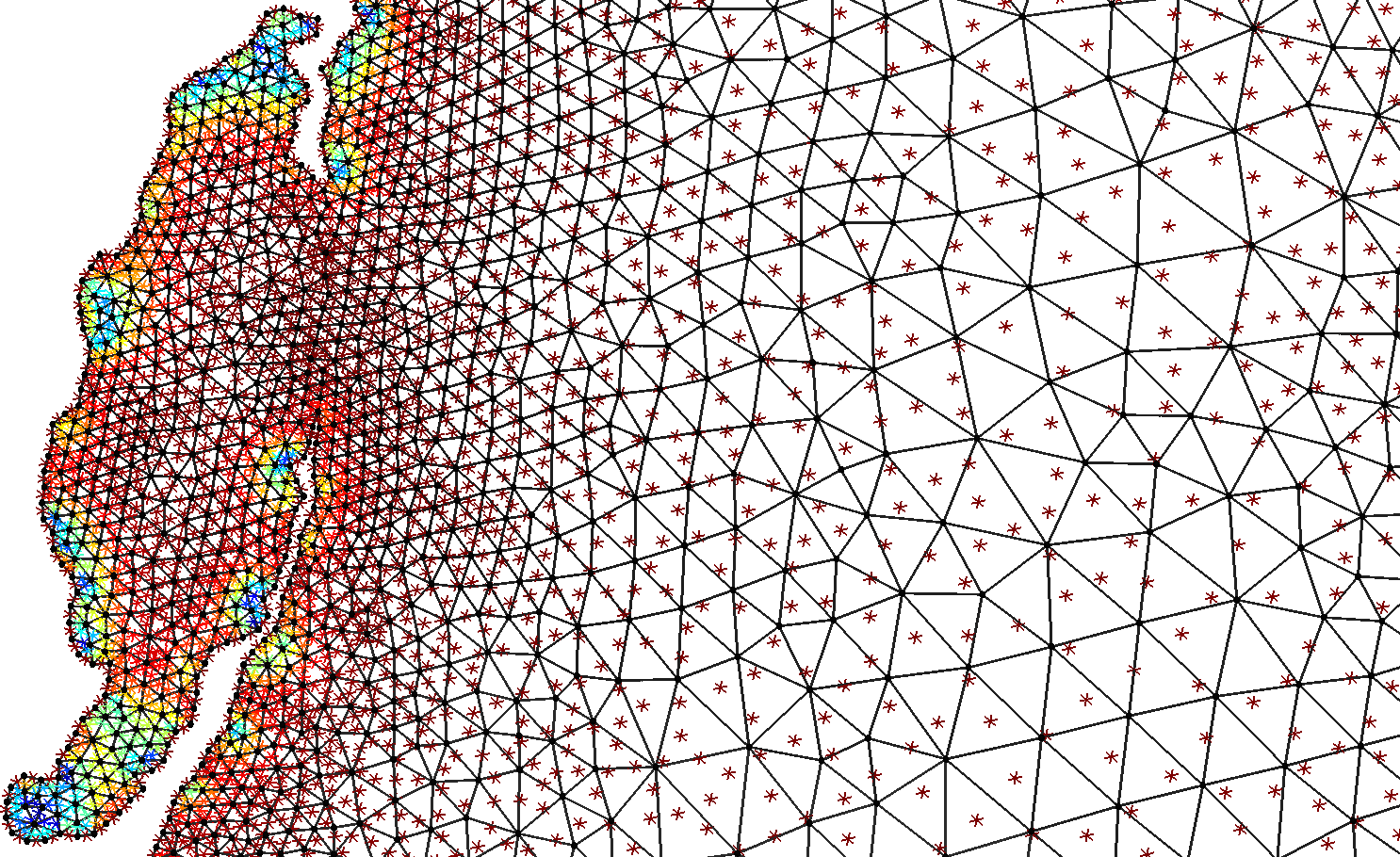}\label{sub:floridaAsterisk}}\\
\subfloat[New triangular mesh.]{\includegraphics[width=.9\linewidth] {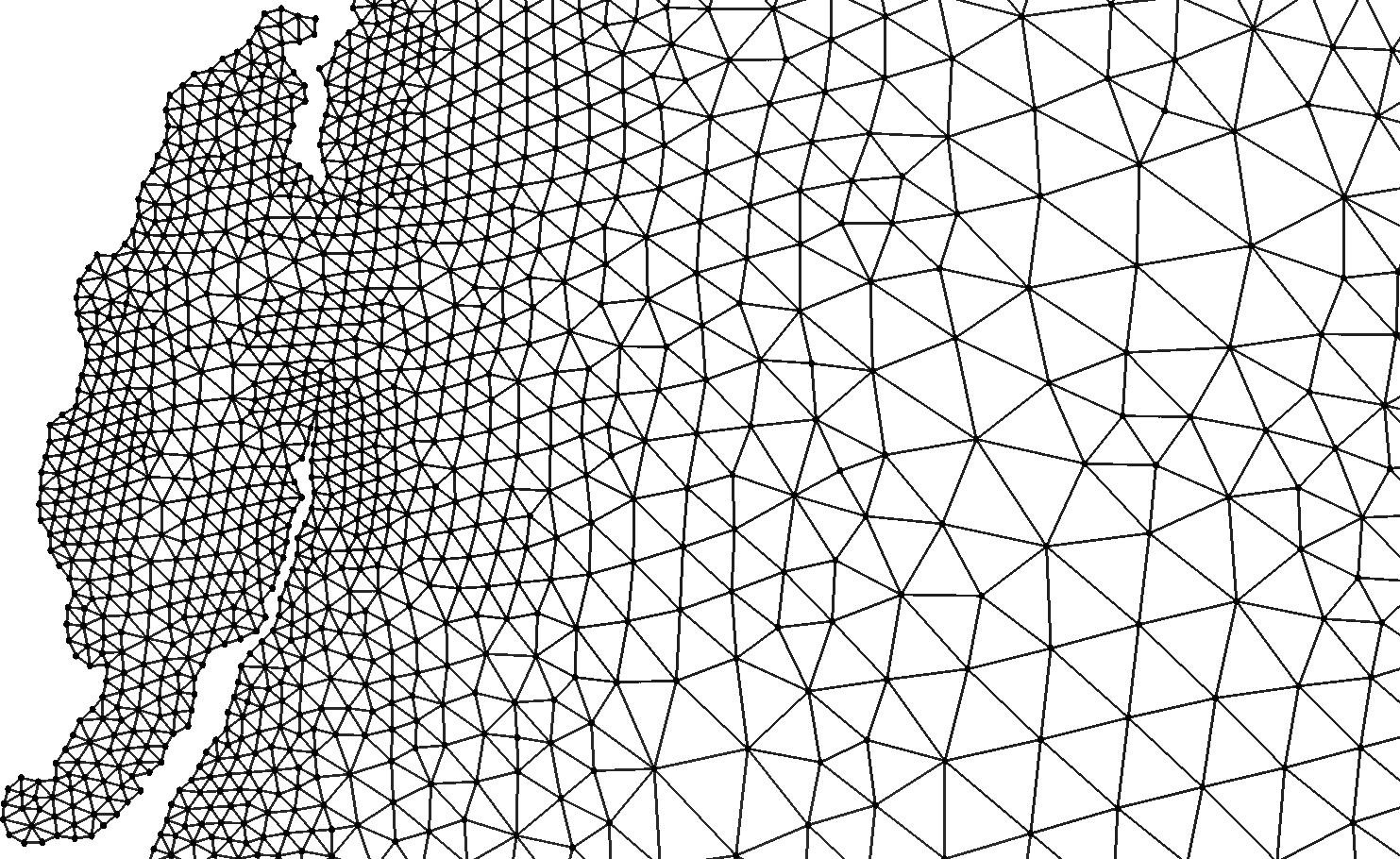}\label{sub:floridaValence}}\\
\end{center}
\caption{Zoom on the Florida keys, the color map is 0 (blue) to 1 (red) and describes the norm of directions. \label{fig:floridaTri}}
\end{figure}

\begin{figure}[!ht]
\begin{center}
\subfloat[Quality value $\gamma$ of initial triangles.]{\includegraphics[width=.5\linewidth] {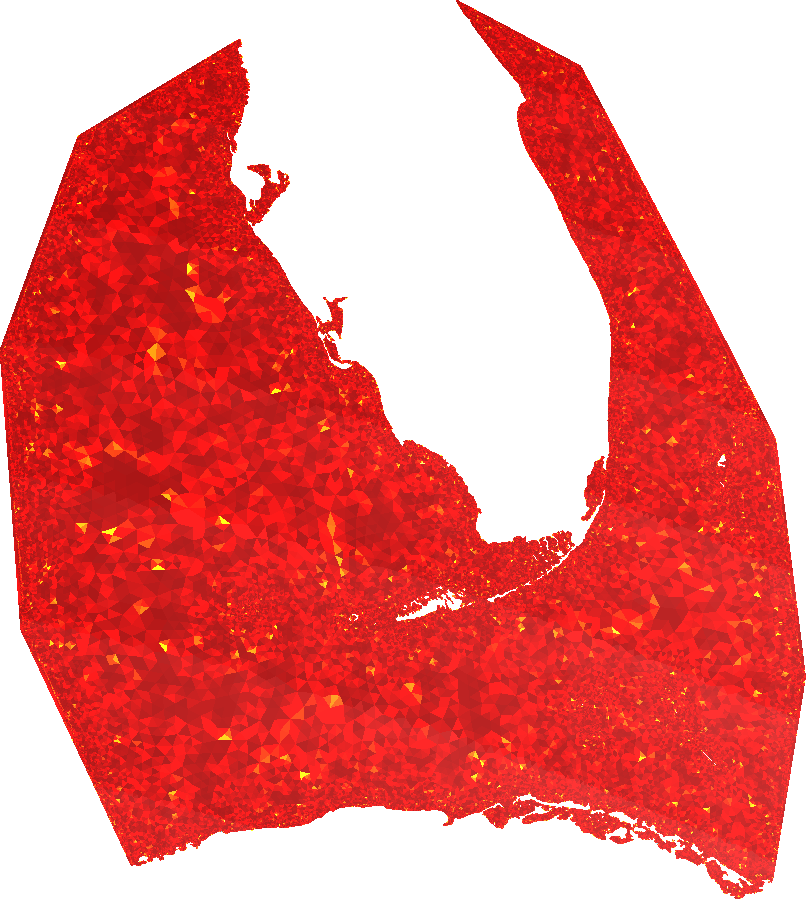}\label{sub:floridaInit}}

\subfloat[Quality value $\gamma$ of final triangles.]{\includegraphics[width=.5\linewidth] {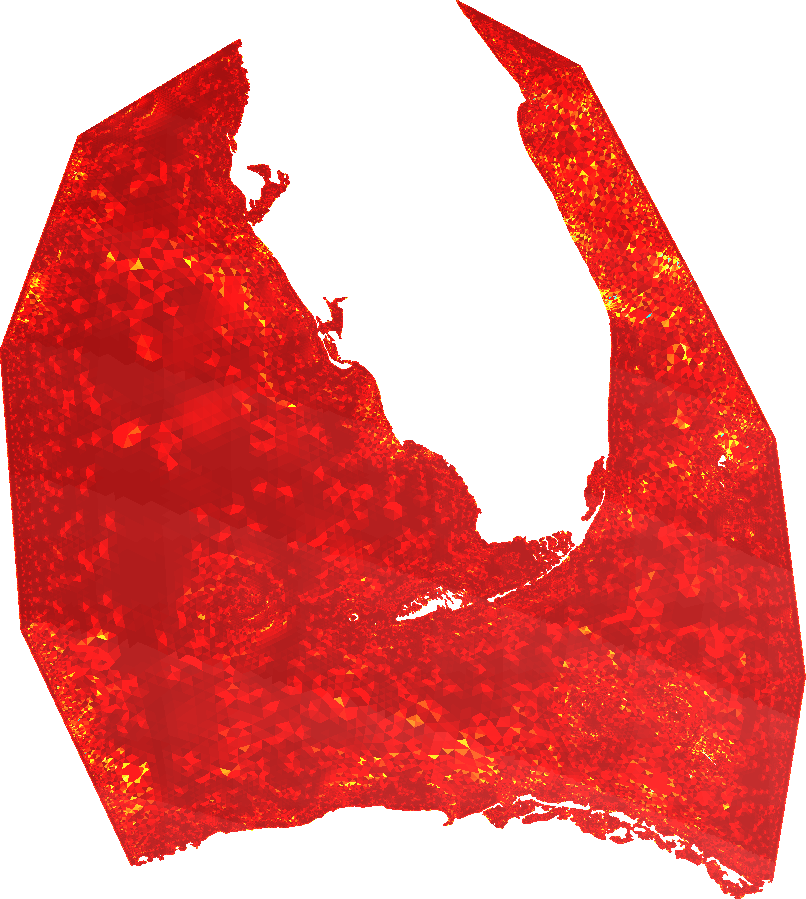}\label{sub:floridaFinal}}

\end{center}
\caption{Zoom on the Florida keys, the color map is 0 (blue) to 1 (red) and describes the quality of triangles. \label{fig:floridaQuality}}
\end{figure}

The input triangular mesh can be improved by using an asterisk field. This field is used to spawn vertices which are consistent with an equilateral triangular grid, Fig. \ref{sub:floridaAsterisk}. The vertices tend to have the correct valence, except in some regions where the size field changes, Fig. \ref{sub:floridaValence}. The final triangular mesh exhibits a smoother distribution of equilateral triangles through the domain, while the mean quality $\bar{\gamma}$ has been improved to 0.9559 (from $\bar{\gamma} = $0.9505 for the initial mesh), Fig. \ref{fig:floridaQuality}. 

Further work will focus on highly improving the numerical scheme that solves Ginzburg-Landau equations, in order to make it competitive.

\FloatBarrier
\section*{Acknowledgements}

The present study was carried out in the framework of the project "Large Scale Simulation of Waves in Complex
Media", which is funded by the Communaut\'e Fran\c caise de Belgique under contract ARC WAVES 15/19-03 with the
aim of developing and using Gmsh~(\cite{geuzaine2009gmsh}).

\bibliography{bibliography}

\begin{thebibliography}{10}
\newcommand{\enquote}[1]{``#1''}
\expandafter\ifx\csname url\endcsname\relax
  \def\url#1{\texttt{#1}}\fi
\expandafter\ifx\csname urlprefix\endcsname\relax\def\urlprefix{URL }\fi

\bibitem{bommes2009mixed}
Bommes D., Zimmer H., Kobbelt L.
\newblock \enquote{Mixed-integer quadrangulation.}
\newblock \emph{ACM Transactions On Graphics (TOG)}, vol.~28, no.~3, 77, 2009

\bibitem{kalberer2007quadcover}
K{\"a}lberer F., Nieser M., Polthier K.
\newblock \enquote{Quadcover-surface parameterization using branched
  coverings.}
\newblock \emph{Computer graphics forum}, vol.~26, pp. 375--384. Wiley Online
  Library, 2007

\bibitem{kowalski2013pde}
Kowalski N., Ledoux F., Frey P.
\newblock \enquote{A PDE based approach to multidomain partitioning and
  quadrilateral meshing.}
\newblock \emph{Proceedings of the 21st international meshing roundtable}, pp.
  137--154. Springer, 2013

\bibitem{cohen2006designing}
Cohen Y.T.P.A.D., Desbrun S.M.
\newblock \enquote{Designing quadrangulations with discrete harmonic forms.}
\newblock \emph{Eurographics symposium on geometry processing}, pp. 1--10. 2006

\bibitem{lee1994new}
Lee C., Lo S.
\newblock \enquote{A new scheme for the generation of a graded quadrilateral
  mesh.}
\newblock \emph{Computers \& structures}, vol.~52, no.~5, 847--857, 1994

\bibitem{remacle2012blossom}
Remacle J.F., Lambrechts J., Seny B., Marchandise E., Johnen A., Geuzainet C.
\newblock \enquote{Blossom-Quad: A non-uniform quadrilateral mesh generator
  using a minimum-cost perfect-matching algorithm.}
\newblock \emph{International journal for numerical methods in engineering},
  vol.~89, no.~9, 1102--1119, 2012

\bibitem{vaxman2016directional}
Vaxman A., Campen M., Diamanti O., Panozzo D., Bommes D., Hildebrandt K.,
  Ben-Chen M.
\newblock \enquote{Directional field synthesis, design, and processing.}
\newblock \emph{Computer Graphics Forum}, vol.~35, pp. 545--572. Wiley Online
  Library, 2016

\bibitem{lai2009metric}
Lai Y.K., Jin M., Xie X., He Y., Palacios J., Zhang E., Hu S.M., Gu X.
\newblock \enquote{Metric-driven rosy field design and remeshing.}
\newblock \emph{IEEE Transactions on Visualization and Computer Graphics},
  vol.~16, no.~1, 95--108, 2009

\bibitem{palacios2007rotational}
Palacios J., Zhang E.
\newblock \enquote{Rotational symmetry field design on surfaces.}
\newblock \emph{ACM Transactions on Graphics (TOG)}, vol.~26, p.~55. ACM, 2007

\bibitem{eppstein2009nineteen}
Eppstein D.
\newblock \enquote{Nineteen proofs of Euler’s formula: V- E+ F= 2.}
\newblock \emph{Information and Computer Sciences, University of California,
  Irvine}, 2009

\bibitem{ray2008n}
Ray N., Vallet B., Li W.C., L{\'e}vy B.
\newblock \enquote{N-symmetry direction field design.}
\newblock \emph{ACM Transactions on Graphics (TOG)}, vol.~27, no.~2, 10, 2008

\bibitem{bethuel1994ginzburg}
Bethuel F., Brezis H., H{\'e}lein F., et~al.
\newblock \emph{Ginzburg-Landau Vortices}, vol.~13.
\newblock Springer, 1994

\bibitem{crouzeix1973conforming}
Crouzeix M., Raviart P.A.
\newblock \enquote{Conforming and nonconforming finite element methods for
  solving the stationary Stokes equations I.}
\newblock \emph{Revue fran{\c{c}}aise d'automatique informatique recherche
  op{\'e}rationnelle. Math{\'e}matique}, vol.~7, no.~R3, 33--75, 1973

\bibitem{ray2016practical}
Ray N., Sokolov D., L{\'e}vy B.
\newblock \enquote{Practical 3d frame field generation.}
\newblock \emph{ACM Transactions on Graphics (TOG)}, vol.~35, no.~6, 233, 2016

\bibitem{saff1997distributing}
Saff E.B., Kuijlaars A.B.
\newblock \enquote{Distributing many points on a sphere.}
\newblock \emph{The mathematical intelligencer}, vol.~19, no.~1, 5--11, 1997

\bibitem{dragnev2002separation}
Dragnev P., Legg D., Townsend D.
\newblock \enquote{On the separation of logarithmic points on the sphere.}
\newblock \emph{Approximation Theory X: Abstract and Classical Analysis}, pp.
  137--144. Vanderbilt University Press, Nashville, TN, 2002

\bibitem{fekete1923verteilung}
Fekete M.
\newblock \enquote{{\"U}ber die Verteilung der Wurzeln bei gewissen
  algebraischen Gleichungen mit ganzzahligen Koeffizienten.}
\newblock \emph{Mathematische Zeitschrift}, vol.~17, no.~1, 228--249, 1923

\bibitem{jezdimirovic2017elliptic}
Jezdimirovic J., Chemin A., Beaufort P., Remacle J.
\newblock \enquote{Elliptic Fekete points obtained by Ginzburg-Landau PDE.}
\newblock Research Notes, 26th International Meshing Roundtable, Sandia
  National Laboratories, September 18-21 2017, 2017

\bibitem{geuzaine2009gmsh}
Geuzaine C., Remacle J.F.
\newblock \enquote{Gmsh: A 3-D finite element mesh generator with built-in
  pre-and post-processing facilities.}
\newblock \emph{International journal for numerical methods in engineering},
  vol.~79, no.~11, 1309--1331, 2009

\end{thebibliography}

\FloatBarrier
\appendix

\section*{Appendix}
This appendix provides an analytical form of the renormalized energy $W(X^c)$
of equation \eqref{eq:phineumann} for a unit disk	. We first compute $\Phi$ and
then prove that second term of \eqref{eqn:renormalized} is equal to
zero in the case of a unit circle. 

\section{Solving the Neumann Problem i.e. Computing $\Phi(x)$ of Equation \eqref{eq:phineumann}}\label{app:analytic}
  
Assume a unit circle $S$. The analytical value of $f$ on the boundary
$\partial S$ of $S$ is  
$$
f = \exp(i~d\theta) ~\text{on} ~ \partial S
$$
as one direction has to be aligned with $\tau$ along the circle.
The Neumann boundary condition is thus
\begin{equation}\label{eqn:neumannbc}
\dfrac{\partial \Phi}{\partial \nu} = d ~ \text{on} ~ \partial \surface
\end{equation}
since $f \times \nabla f \cdot \tau = d$ on $\partial \surface$.
Indeed, from 
$$
a\bar{b} = a \cdot b - i ~ a \times b, \forall a,b \in \mathbb{C}
$$
and
$$
\nabla f \cdot \tau = d~i~f \nabla \theta \cdot \tau = i ~df
$$
the condition (\ref{eqn:neumannbc}) corresponds to the imaginary part of the corresponding complex product.

The Green function of the two-dimensional Laplacian operator is part of $\Phi$
\begin{equation}
\Phi^0 (\mathbf{x}) = \sum_{i=1}^d \log|\mathbf{x} - \mathbf{x}^c_i|
\end{equation}
Even if $\nabla^2 \Phi^0 = 2\pi \sum_{i=1}^d \delta (\mathbf{x} - \mathbf{x}_i^c)	$, the flux (per unit of length) $\nabla \Phi^0 \cdot \nu$  does not correspond to (\ref{eqn:neumannbc}).
The solution $\Phi$ contains another term $\Phi^1$.
It may also be written as a sum of the contributions coming from the $d$ critical points.
Therefore,
$$
\Phi = \Phi^0 + \Phi^1 = \sum_{i=1}^d \left( \Phi^0_i + \Phi^1_i \right)
$$
such that
\begin{equation}\label{eqn:phi1}
\left.\begin{array}{rcl}
\Phi^0_i(\mathbf{x}) &=& \log|\mathbf{x} - \mathbf{x}^c_i| \\
\nabla^2 \Phi^1_i &=& 0 ~ \text{in}  ~  \surface \\
(\nabla \Phi^0_i + \nabla \Phi^1_i) \cdot \nu &=& 1 ~  \text{on} ~ \partial \surface
\end{array}\right.
\end{equation}
Function $\Phi^1_i$ can be written as series of circular harmonics 
$$
\Phi^1_i(r;\theta) = A_{i,0} + \sum_{n=1}^\infty r^n \left[ A_{i,n} \cos(n~\theta) + B_{i,n} \sin(n~\theta)  \right]
$$
where $(r;\theta)$ are polar coordinates.
We search for the solution of a Neumann
problem which is defined to a constant. Setting set $A_{i,0}=0$
assigns to zero the average of $\Phi^1_i$.
The idea is simple. We use $\Phi^1$ which is harmonic to remove all
oscillatory parts of $\nabla \Phi \cdot \nu$ 
along the boundary $\partial \surface$.

Let us assume the i-th critical point is located on the x axis (the
real axis), i.e. $\mathbf{x}_i^c=(r^c;0)$ with 
cartesian coordinates (see Figure \ref{fig:unitdisk}).
\begin{figure}
\begin{center}
\begin{tikzpicture}
\draw[thick] (0,0) circle (2cm);
\draw[very thick,->] (-.5,0)--(2.5,0);
\draw[very thick,->] (0,-.5)--(0,2.5);
\draw[thick,red] (0,0)--(1.7,0);
\draw[thick,blue] (0:0)--(70:1cm);
\draw[thick,green] (0:1.7cm)--(70:1cm);
\draw[thick,->] (.5,0) arc (0:70:.5cm);
\node at (1.7,0) {$\bullet$};
\node at (1.7,.25) {$\mathbf{x}_i^c$};
\node at (70:1cm) {$\bullet$};
\node at (73:1.25cm) {$\mathbf{x}$};
\node at (35:.75cm) {$\theta$};
\node[red] at (.85,-.25) {$r^c$};
\node[blue] at (80:.5) {$r$};
\end{tikzpicture}
\end{center}
\caption{Unit disk where the i-th critical point is depicted.}\label{fig:unitdisk}
\end{figure}
One has 
$$|\mathbf{x}-\mathbf{x}_i^c|^2 = r^2 + {r^c}^2 - 2 r r^c \cos (\theta)$$
and
$$
\left. {\partial \Phi^0_i  \over \partial r} \right|_{r=1}
v= \left. {1 \over |\mathbf{x} - \mathbf{x}_i^c|} {\partial |\mathbf{x} - \mathbf{x}_i^c| \over \partial r}\right|_{r=1} =
{1 - r^c \cos \theta \over 1 + {r^c}^2 - 2 r^c \cos \theta}. 
$$
The last expression can be reformulated as
$$ 
{(1 + {r^c}^2 - 2 r^c \cos \theta )+ 1-{r^c}^2 \over 2( 1 + {r^c}^2 - 2 r^c \cos \theta)}
= {1 - s \cos \theta + {1-{r^c}^2 \over 1+{r^c}^2} \over 2( 1 -s \cos \theta)}
$$
with $s = 2r^c/(1 + {r^c}^2)$.
Taking into account the identity 
$$ { 1 \over 1 -x} = \sum_{n=0}^{\infty} x^n, ~~~|x| < 1,$$
we have
\begin{equation}
\left. {\partial \Phi^0_i  \over \partial r} \right|_{r=1}
= {1 \over 1 + {r^c}^2} \left(1 
+ \left[\sum_{n=1}^\infty s^n (\cos \theta)^n \right]{1-{r^c}^2 \over 2} \right)
\label{eq:pow}
\end{equation}
and
\begin{equation}
\left. {\partial \Phi_i^1  \over \partial r} \right|_{r=1}  =
\sum_{n=1}^{\infty} n (A_{i,n} \cos (n\theta) + B_{i,n} \sin (n\theta)).
\label{eq:n}
\end{equation}

Powers of $\cos (\theta)$ appear in \eqref{eq:pow}. In order to
replace such powers by $\cos( n ~\theta)$'s like in Equation \eqref{eq:n}, 
we use a well known property of Chebyshev polynomials\,:
$P_n(\cos( \theta)) = \cos (n ~\theta)$.
We thus have
\begin{equation}
\cos (m \theta)  = \sum_{n=0}^m {\mathcal P}_{mn} (\cos \theta)^n
\label{eq:cheb}
\end{equation}
where the  ${\mathcal P}_{mn}$'s are the entries 
of the Chebyshev coefficient matrix $\mathcal P$. 
Equation \eqref{eq:cheb} can thus be regarded as a system of equations
$$
X_i = \mathcal{P}_{in} Y_n
\quad , \quad 
X_i = \cos (i\theta)
\quad , \quad 
Y_n = (\cos\theta)^n.
$$
The system matrix $\mathcal P$ is lower triangular, 
so the system can be inverted easily
$$
Y_n =  \mathcal{P}^{-1}_{ni} X_i,
$$
or equivalently, back with the initial notation, 
$$(
\cos (\theta))^n  = \sum_{i=0}^n {\mathcal P}^{-1}_{ni} (\cos (i\theta)).
$$ 
Thus, 
$$\sum_{n=1}^\infty s^n (\cos \theta)^n = \sum_{n=1}^\infty s^n \sum_{i=0}^n {\mathcal P}^{-1}_{ni} (\cos i \theta) = \sum_{n=0}^\infty w_n \cos n \theta$$
with
$$w_n = \sum_{i=n}^\infty s^i  {\mathcal P}^{-1}_{in}.$$
Finally, we get the following series for the normal derivative of $\Phi_i^0$:
$$ 
\left. {\partial \Phi_i^0  \over \partial r} \right|_{r=1}
=  \underbrace{{1 \over 1 + {r^c}^2} \left(1 +  w_0{1-{r^c}^2 \over 2} \right)}_{W_0=1} +
\sum_{n=1}^\infty \underbrace{{1-{r^c}^2 \over 2(1 + {r^c}^2)}  w_n}_{W_n} (\cos n\theta)
$$ 
We get the final condition
\begin{eqnarray}
\left. {\partial (\Phi_i^0+\Phi_i^1)  \over \partial r} \right|_{r=1}
=
W_0 + \sum_{n=1}^\infty \left[(W_n + n A_{i,n})\cos (n\theta) + n B_{i,n} \sin(n \theta)\right].
\end{eqnarray}
The boundary condition should be non oscillatory:
So, 
$$B_{i,n}=0~~~ \mbox{and}~~~A_{i,n} = -{W_n \over n}.$$
Finally
$$
\Phi(\mathbf{x}) = \sum_{i=1}^d 
\left[\log |\mathbf{x} - \mathbf{x}_i^c|
+   \sum_{n=1}^{\infty} r^n A_{i,n} \cos (n\theta)
\right] .  
$$

\section{$\displaystyle{\int_{\partial \surface} \Phi~ f \times f_\tau ~ds}$ is Zero along a Circle}\label{app:integral}

We want to show that
$$
\int_{\partial \surface} \Phi ~ f \times \nabla f \cdot \tau ~ ds = 0
$$
when $\partial \surface$ is a unit circle.

In the previous section, we have shown that
$$
f \times \nabla f \cdot \tau = d
$$
Besides, $\Phi^1$ has been derived such that it is non oscillatory along the unit circle.
Hence, it remains to show
$$
\sum_{k=1}^d \int_{\partial \surface} \log|\mathbf{x}-\mathbf{x}_k^c| ~ ds = 0
$$

We can express that integral with complex variables
$$
\Re \left\{ \sum_{k=1}^d \oint_{|z|=1} \log(z-z_k^c) ~ dz  \right\} = 0
$$
with the complex logarithm, which has two features:
\begin{itemize}
\item the complex logarithm is a multivalued function
\item the complex logarithm has a peculiar singularity in zero
\end{itemize}

Those features are due to the fact that zero is a branch point.
In our case, the branch points are the critical points $z_k^c$.
A branch cut has to be drawn for each critical point.
If $z_k^c = r^c \exp(i~\gamma_k)$, the branch cut is such that $\theta \in [\gamma_k;\gamma_k+2\pi), z=r\exp(i~\theta) \in \mathbb{C}$ (red line on Fig. \ref{fig:log}).
The logarithm is analytical thanks to the branch cut, and its contour integral may be written as
$$
{\scriptsize
\begin{array}{rcl}
  \displaystyle{\oint_{|z|=1} \log\left(z-z_k^c\right) ~ dz}  &=&  \displaystyle{{\color{white}+} \int_0^{\gamma_k^-} \log\left(\exp(i~t) - r^c \exp(i~\gamma_k)\right) i \exp(i~t) ~ dt} \\
&& \displaystyle{+ \int_{\gamma_k^+}^{2\pi} \log\left(\exp(i~t) - r^c \exp(i~\gamma_k)\right) i \exp(i~t) ~ dt}
\end{array}}
$$

The antiderivative of $\log(w)$ is $w \left(\log(w) - 1\right)$, yielding the following
$$
{\fontsize{5}{6}\selectfont 
\begin{array}{rcl}
\displaystyle{\oint_{|z|=1} \log\left(z-z_k^c\right) ~ dz} &=& \displaystyle{{\color{white}+}\left[ \left(\exp(i~t) - r^c \exp(i~\gamma_k)\right) \left(\log\left(\exp(i~t) - r^c \exp(i~\gamma_k)\right) - 1\right) \right]_0^{\gamma_k^-}}  \\
 && \displaystyle{+ \left[ \left(\exp(i~t) - r^c \exp(i~\gamma_k)\right) \left(\log\left(\exp(i~t) - r^c \exp(i~\gamma_k)\right) - 1\right) \right]_{\gamma_k^+}^{2\pi}} \\
 &=& \displaystyle{\exp(i ~\gamma_k) (1-r^c) \left( \log\left(\exp(i~\gamma_k^-)\right) - \log\left(\exp(i~\gamma_k^+)\right)  \right)}\\
 &=& 2 \pi i ~ \exp(i ~\gamma_k) (1-r^c)
\end{array}}
$$
The last step is due to the fact that $\gamma_k^- - \gamma_k^+ = 2\pi$ (the branch cut, in other words).

$$
\oint_\mathcal{C} \log(w) ~ dw = 2\pi i ~ \exp(i~\gamma_k) (1-r^c)
$$
which depends on the branch cut.

\begin{figure}
\begin{center}
\begin{tikzpicture}[scale=.9]
\draw[thick,dotted] (0,0) circle (4cm);
\draw[thick,->] (-.5,0)--(4.5,0);
\draw[thick,->] (0,-.5)--(0,4.5);
\draw[very thick,red] (55:3)--(55:8);
\node at (60:3.15) {$z_k^c$};
\node at (55:3) {$\bullet$};
\draw[thick,->,ForestGreen] (3.25,0) arc (0:55:3.25);
\node[ForestGreen] at (27.5:3.5) {$\gamma_k$};
\end{tikzpicture}
\end{center}
\caption{Illustration of the k-th critical point defining a branch cut (in red) for the corresponding complex logarithm.}\label{fig:log}
\end{figure}

The $d$ critical points $z_k^c$ being at the same distance $r^c$ from the origin and evenly spaced by an angles $2\pi/d$, we get
$$
\sum_{k=1}^d \oint_{|z|=1} \log(z-z_k^c) ~  dz = 2\pi i (1-r^c) ~ \sum_{k=1}^d \exp(i~[\gamma_1 + (k-1)2\pi/d])
$$
which is zero since the sum of $d>1$ complex numbers corresponding to points evenly distributed along on a circle centered at the origin is zero.

\end{document}